\documentclass[11pt]{article}
\usepackage{t1enc}
    \usepackage[latin1]{inputenc}
    \usepackage[english]{babel}
        \usepackage{latexsym}
\begin{document}
\setlength{\arraycolsep}{.136889em}
\renewcommand{\theequation}{\thesection.\arabic{equation}}
\newtheorem{thm}{Theorem}[section]
\newtheorem{propo}{Proposition}[section]
\newtheorem{lemma}{Lemma}[section]
\newtheorem{corollary}{Corollary}[section]
\centerline{\Large\bf Random walk local time approximated by a Wiener
sheet}
\medskip
\centerline{\Large\bf combined with an independent Brownian motion}

\bigskip\bigskip
\bigskip\bigskip
\renewcommand{\thefootnote}{1}
\noindent
{\textbf{Endre Cs\'{a}ki}\footnote{Research supported by the
Hungarian National Foundation for Scientif\/ic Research, Grant No.
K 61052 and K 67961.}}
\newline
Alfr\'ed R\'enyi Institute of
Mathematics, Hungarian Academy of Sciences, Budapest, P.O.B. 127,
H-1364, Hungary. E-mail address: csaki@renyi.hu

\bigskip
\renewcommand{\thefootnote}{2}
\noindent
{\textbf{Mikl\'os Cs\"org\H{o}}
\footnote{Research supported by an NSERC Canada Discovery Grant at
Carleton University}}
\newline
School of Mathematics and Statistics, Carleton University,
1125 Colonel By Drive, Ottawa, Ontario, Canada K1S 5B6.
 E-mail address: mcsorgo@math.carleton.ca

\bigskip
\renewcommand{\thefootnote}{3}
\noindent
{\textbf{Ant\'{o}nia
F\"{o}ldes}\footnote{Research supported by a PSC CUNY Grant, No.
68030-0037.}}
\newline
Department of Mathematics, College of Staten
Island, CUNY, 2800 Victory Blvd., Staten Island, New York 10314,
U.S.A.  E-mail address: foldes@mail.csi.cuny.edu

\bigskip
\noindent
{\textbf{P\'al R\'ev\'esz}$^1$}
\newline
Institut f\"ur
Statistik und Wahrscheinlichkeitstheorie, Technische Universit\"at
Wien, Wiedner Hauptstrasse 8-10/107 A-1040 Vienna, Austria.
E-mail address: reveszp@renyi.hu

\bigskip\centerline
{\it Dedicated to the memory of Walter Philipp}

\medskip
\begin{abstract}
Let $\xi(k,n)$ be the local time of a simple symmetric random walk on
the line. We give a strong approximation of the centered local time
process $\xi(k,n)-\xi(0,n)$ in terms of a Wiener sheet and an independent
Wiener process, time changed by an independent Brownian local time. Some
related results and consequences are also established.
\end{abstract}


\bigskip

 \vspace{.1cm}

\section{Introduction and main results}
\renewcommand{\thesection}{\arabic{section}} \setcounter{equation}{0}
\setcounter{thm}{0} \setcounter{lemma}{0}

Let $X_i$, $i=1,2,\ldots$, be i.i.d. random variables with the
distribution $P(X_i=1)=P(X_i=-1)=1/2$ and put $S_0:=0,$
$S_i:=X_1+\ldots+X_i,\, \, i=1,2,\ldots$. Define the local time process
of this simple symmetric random walk by
\begin{equation}
\xi(k,n):=\#\{i:\, 1\leq i\leq n,\, S_i=k\},\quad k=0,\pm 1,\pm
2,\ldots,
\, n=1,2,\ldots
\label{loc}
\end{equation}
We can also interpret $\xi(k,n)$ as the number of excursions
away from $k$ completed before $n$.

We now define some related quantities for further use. Let $\rho_0:=0$
and
\begin{equation}
\rho_i:=\min\{j>\rho_{i-1}:\, S_j=0\},\quad i=1,2,\ldots,
\label{rho}
\end{equation}
i.e., $\rho_i$ is the time of the $i$-th return to zero, or, in
other words, the endpoint of the $i$-th excursion away from $0$.
We say that $(S_a,S_{a+1},\ldots,S_b)$ is an excursion away from $k$,
if $S_a=S_b=k,\, S_i\neq k,\, a<i<b$. This excursion will be called
upward if $S_i>k,\, a<i<b$ and downward if $S_i<k,\, a<i<b$.
Define $\rho_i^+$ as the endpoint of the $i$-th upward excursion away
from $0$, and let $\xi(k,n,\uparrow)$ be the number of upward excursions
away from $k$ completed up to time $n$. Similarly, let
$\xi(k,n,\downarrow)$ be the number of downward excursions away from $k$
completed up to time $n$.

Let $\{W(t),\, t\geq 0\}$ be a standard Wiener process and consider its
two-parameter local time process $\{\eta(x,t),\, x\in R,\, t\geq 0\}$
satisfying
\begin{equation}
\int_A\eta(x,t)\, dx=\lambda\{s:\, 0\leq s\leq t,\, W(s)\in A\}
\label{eta}
\end{equation}
for any $t\geq 0$ and Borel set $A\subset R$, where $\lambda(\cdot)$ is
the Lebesgue measure. In the sequel we simply call $\eta(\cdot,\cdot)$
a standard Brownian local time.

The study of the asymptotic behaviour of the centered local time
processes $\xi(k,n)-\xi(0,n)$ and $\eta(x,t)-\eta(0,t)$ has played a
significant role in the development of the local time theory of
random walks and that of Brownian and iterated Brownian motions. The
first of this kind of results we have in mind is due to Dobrushin
\cite{do}. Namely, in this landmark paper, a special case of one of his
theorems for additive functionals of a simple symmetric random walk
reads as follows.

\medskip\noindent
\textbf{Theorem A1} {\it For any $k=1,2,\ldots$}
\begin{equation}
\frac{\xi(k,n)-\xi(0,n)}{(4k-2)^{1/2}n^{1/4}}\to_d
U\sqrt{|V|}, \quad n\to\infty,
\label{dob}
\end{equation}
{\it where $U$ and $V$ are two independent standard normal variables}.

Here and in the sequel $\to_d$ denotes convergence in
distribution.

\medskip
On the other hand, concerning now centered Brownian local times, a
special case of a more general fundamental theorem of Skorokhod and
Slobodenyuk \cite{sks1}, that is an analogue of Dobrushin's theorem as
in \cite{do}, yields the following result.

\medskip\noindent
\textbf{Theorem B1} {\it For any $x> 0$}
\begin{equation}
\frac{\eta(x,t)-\eta(0,t)}{2x^{1/2}t^{1/4}}\to_d
U\sqrt{|V|}, \quad t\to\infty,
\label{skos}
\end{equation}
{\it where $U$ and $V$ are two independent standard normal variables}.

\medskip
While these two theorems are similar, we call attention to their
intriguing difference in their scaling constants. For example, the
respective scaling constant for $k=x=1$ is $2^{1/2}$ in Theorem A1 and
it is $2$ in Theorem B1.

Dobrushin's result as in \cite{do} was extended under various conditions
by Kesten \cite{ke}, Skorokhod and Slobodenyuk \cite{sks2}, Kasahara
\cite{ka1}, \cite{ka2} and Bo\-ro\-din \cite{bo2}. For some details on
the nature of these extensions we refer to the Introduction in \cite{cscsfr2}.

In connection with the analogue of (\ref{dob}) as spelled out in
(\ref{skos}), for further extensions along these lines we refer to
Papanicolaou et al. \cite{psv}, Ikeda and Watanabe \cite{iw} and the
survey paper of Borodin \cite{bo}. For some details we again refer to
\cite{cscsfr2}.

The papers mentioned in the previous two paragraphs, in general, are
concerned with studying additive functionals of the form
$A_n:=\sum_{i=1}^n f(S_i)$, and their integral forms
$I_t:=\int_0^t g(W(s))\, ds$, where
$f(x),\, x\in R$, and $g(x),\, x\in R$, are real valued functions
satisfying appropriate conditions. In particular, Cs\'aki et al.
\cite{cscsfr2} deals with strong approximations of these two types of
additive functionals, together with their weak and strong convergence
implications.

In view of (\ref{dob}) and (\ref{skos}) above, we now mention some
corresponding iterated logarithm laws. For example, (4.1a) of Cs\'aki et
al. \cite{cscsfr2} yields

\medskip\noindent
\textbf{Theorem C1} \textit{For} $k=1,2,\ldots$ \textit{we have}
\begin{equation}
\limsup_{n\to\infty}\frac{\xi(k,n)-\xi(0,n)}
{(4k-2)^{1/2}n^{1/4}(\log\log n)^{3/4}}=\frac{2}{3}6^{1/4}\quad
\mathrm{a.s.}
\label{xi2}
\end{equation}

\medskip
While studying the local time process of a symmetric random walk
standardized by its local time at zero, Cs\"org\H o
and R\'ev\'esz \cite{csr85} established the next result.

\medskip\noindent
\textbf{Theorem C2} \textit{For} $k=1,2,\ldots$ \textit{we have}
\begin{equation}
\limsup_{n\to\infty}\frac{\xi(k,n)-\xi(0,n)}{(4k-2)^{1/2}(\xi(0,n)\log\log
n)^{1/2}}=2^{1/2}\quad \mathrm{a.s.}
\label{xi1}
\end{equation}

\medskip
Moreover, Theorem 1 of Cs\'aki and F\"oldes \cite{csf} yields the
next pair of Theorems.

\medskip\noindent
\textbf{Theorem D1} \textit{For} $x>0$ \textit{we have}
\begin{equation}
\limsup_{t\to\infty}\frac{\eta(x,t)-\eta(0,t)}
{2x^{1/2}t^{1/4}(\log\log t)^{3/4}}=\frac{2}{3}6^{1/4}
\quad\mathrm{a.s.}
\label{eta2}
\end{equation}

\medskip\noindent
and

\noindent
\textbf{Theorem D2} \textit{For} $x>0$ \textit{we have}
\begin{equation}
\limsup_{t\to\infty}\frac{\eta(x,t)-\eta(0,t)}
{2x^{1/2}(\eta(0,t)\log\log t)^{1/2}}=2^{1/2} \quad\mathrm{a.s.}
\label{eta1}
\end{equation}

\medskip
While these two pairs of theorems are similar, just like in case of
(\ref{dob}) and (\ref{skos}), we call attention to their intriguing
difference in their scaling constants.

In view of Theorems C2 and D2, we state the next two results.

\medskip\noindent
\textbf{Theorem A2} \textit{For} $k=1,2,\ldots$ \textit{we have}
\begin{equation}
\frac{\xi(k,n)-\xi(0,n)}{(4k-2)^{1/2}(\xi(0,n))^{1/2}}\to_d
U, \quad n\to\infty,
\label{xi3}
\end{equation}
\textit{where $U$ is a standard normal random variable.}

\medskip\noindent
\textbf{Theorem B2} \textit{For} $x>0$ \textit{we have}
\begin{equation}
\frac{\eta(x,t)-\eta(0,t)}{2x^{1/2}(\eta(0,t))^{1/2}}\to_d
U, \quad t\to\infty,
\label{eta3}
\end{equation}
\textit{where $U$ is a standard normal random variable}.

\medskip
Theorem A2 is argued \textit{intuitively} on p. 90 of Cs\"org\H o and
R\'ev\'esz \cite{csr85}, and it can be rigorously based on our results
in Cs\'aki et al. \cite{cscsfr2}, while Theorem
B2 is stated as one of the consequences of our results in Cs\'aki et al.
\cite{cscsfr1}.

The next weak convergence result for fixed $k$ follows from
Kasahara \cite{ka1}.

\medskip\noindent
\textbf{Theorem E} \textit{For $k=1,2,\ldots$ we have}
$$
\frac{\xi(k,[\lambda t])-\xi(0,[\lambda t])}
{(4k-2)^{1/2}\lambda^{1/4}}
\to_w W(\widetilde\eta(0,t)),\qquad \lambda\to\infty,
$$
\textit{where $\widetilde\eta(\cdot,\cdot)$ is a standard Brownian local
time, independent of the Wiener process $W(\cdot)$.}

Here and in the sequel $\to_w$ denotes weak convergence in the respective
function spaces in hand (here $D[0,\infty)$).

\medskip
Moreover, for fixed $x$ the next weak convergence result in
$C[0,\infty)$ follows from Papanicolaou et al. \cite{psv}.

\medskip\noindent
\textbf{Theorem F} \textit{For $x>0$ we have}
$$
\frac{\eta(x,\lambda t)-\eta(0,\lambda t)}{2x^{1/2}\lambda^{1/4}}
\to_w W(\widetilde\eta(0,t)),\qquad \lambda\to\infty,
$$
\textit{where $\widetilde\eta(0,t)$ is as in} Theorem E.

When our paper \cite{cscsfr2} on strong approximations of additive
functionals is interpreted in our present context, its general results
also imply strong approximations for $\xi(k,n)-\xi(0,n)$ when $k$ is fixed,
as spelled out in the next theorem.

\medskip\noindent
\textbf{Theorem G} \textit{On an appropriate probability space
for a simple symmetric random walk $\{S_i,\, i=0,1,\ldots\}$, for any
$k=1,2,\ldots$, we can construct a standard Wiener process
$\{W(t),\, t\geq 0\}$ and, independently of the latter, a standard
Brownian local time $\{\widetilde \eta(0,t),\, t\geq 0\}$ such that,
as $n\to\infty$, with sufficiently small $\varepsilon>0$ we have
\begin{equation}
\xi(k,n)-\xi(0,n)=(4k-2)^{1/2}W(\widetilde
\eta(0,n))+O(n^{1/4-\varepsilon})\quad\mathrm{a.s.}
\label{inv1}
\end{equation}
and
\begin{equation}
\xi(0,n)-\widetilde\eta(0,n)=O(n^{1/2-\varepsilon})\quad\mathrm{a.s.}
\label{inv2}
\end{equation}}

\medskip
Following the method of proof of Theorem 2 in Section 3 of
\cite{cscsfr2}, one can also establish the next theorem, which is also a
consequence of our Theorem in \cite{cscsfr1}, that is quoted below (cf.
Theorem J).

\medskip\noindent
\textbf{Theorem H} \textit{On an appropriate probability space
for the standard Brownian local time process $\{\eta(x,t),\, x\in R,\,
t\geq 0\}$ of a standard Brownian motion,  for any
$x>0$, we can construct a standard Wiener process
$\{W(t),\, t\geq 0\}$ and, independently of the latter, a standard
Brownian local time $\{\widetilde \eta(0,t),\, t\geq 0\}$ such that,
as $t\to\infty$, with sufficiently small $\varepsilon>0$ we have
\begin{equation}
\eta(x,t)-\eta(0,t)=2x^{1/2}W(\widetilde
\eta(0,t))+O(t^{1/4-\varepsilon})\quad\mathrm{a.s.}
\label{inv3}
\end{equation}
and
\begin{equation}
\eta(0,t)-\widetilde\eta(0,t)=O(t^{1/2-\varepsilon})\quad\mathrm{a.s.}
\label{inv4}
\end{equation}}

\medskip
A common property of the above quoted theorems is that they treat the
two-time parameter processes $\xi(k,n)$ and $\eta(x,t)$ for $k$,
respectively $x$, fixed, i.e., as if they were one-time parameter
stochastic processes. (In (\ref{inv1}), resp. (\ref{inv3}), both $W$ and
the $O$ term may depend on $k$, resp. $x$.) Clearly, studying them as
two-time parameter processes is of cardinal interest. A significant first
step along these lines was made by Yor \cite{yo}, who established the
following weak convergence result.

\medskip\noindent
\textbf{Theorem I} \textit{As $\lambda\to\infty$},
$$
\left(\frac{1}{\lambda}W(\lambda^2 t),
\frac{1}{\lambda}\eta(x,\lambda^2 t),
\frac{1}{2\sqrt{\lambda}}(\eta(x,\lambda^2 t)-\eta(0,\lambda^2
t))\right)
$$
$$
\to_{w} (W(t), \eta(x,t), W^*(x,\eta(0,t))),
$$
\textit{where $W^*(\cdot,\cdot)$ is a Wiener sheet, independent of the
standard Wiener process  $W(\cdot)$, $\eta(\cdot,\cdot)$ is
the local time of $W(\cdot)$, and $\to_{w}$ denotes weak convergence over
the space of all continuous functions from $R_+^2$ to $R^3$, endowed
with the topology of compact uniform convergence.}

\medskip
By a Wiener sheet we mean a two-parameter Gaussian process
$$\{W(x,y),\, x\geq 0,\, y\geq 0\}$$ with mean $0$ and covariance
function
$$
EW(x_1,y_1)W(x_2,y_2)=(x_1\wedge x_2)(y_1\wedge y_2)
$$
(cf., e.g., Section 1.11 in Cs\"org\H o and R\'ev\'esz \cite{csr}).

In Cs\'aki et al. \cite{cscsfr1} we proved the following strong
approximation of Brownian local time by a Wiener sheet.

\medskip\noindent
\textbf{Theorem J} \textit{On an appropriate probability space
for the standard Brownian local time process $\{\eta(x,t),\, x\in R,\,
t\geq 0\}$ of a standard Brownian motion, we can construct a Wiener
sheet $\{W(x,u),$ $x,u\geq 0\}$ and, independently of the latter, a
standard Brownian local time $\{\widetilde \eta(0,t),\, t\geq 0\}$ such
that, as $t\to\infty$, for sufficiently small $\varepsilon>0$ there
exists $\delta>0$ for which we have
$$
\sup_{0\leq x\leq t^{\delta}}
|\eta(x,t)-\eta(0,t)-2W(x,\widetilde\eta(0,t))|=O(t^{1/4-\varepsilon})
\quad\mathrm{a.s.}
$$
and}
$$
\eta(0,t)-\widetilde\eta(0,t)=O(t^{1/2-\varepsilon})\quad\mathrm{a.s.}
$$

\medskip
In Cs\'aki et al. \cite{cscsfr1} we also proved the analogue of Theorem
J for $t$ replaced by the inverse local time $\alpha(\cdot)$ defined by
$$
\alpha(u):=\inf\{t\geq 0: \, \eta(0,t)\geq u\}.
$$

\medskip\noindent
\textbf{Proposition A} \textit{On an appropriate probability space
for the standard Brownian local time process $\eta(x,t),\, x\in R,\,
t\geq 0\}$ of a standard Brownian motion, we can construct a Wiener
sheet $\{W(x,u),\, x,u\geq 0\}$ and, independently of the latter, an
inverse local time process $\{\widetilde \alpha(u),\, u\geq 0\}$ such
that for sufficiently small $\varepsilon>0$ there exists $\delta>0$ for
which as $u\to\infty$, we have
$$
\sup_{0\leq x\leq u^{\delta}}
|\eta(x,\alpha(u))-u-2W(x,u)|=O(u^{1/2-\varepsilon})
\quad\mathrm{a.s.}
$$
and}
$$
\alpha(u)-\widetilde\alpha(u)=O(u^{2-\varepsilon})\quad\mathrm{a.s.}
$$

Concerning weak convergence of increments of random walk local time,
Eisenbaum \cite{ei} established a two-parameter result for symmetric
Markov chains at inverse local times, which for a simple symmetric
random walk reads as follows.

\medskip\noindent
\textbf{Proposition B} \textit{As $\lambda\to\infty$,
$$
\frac{\xi(k,\rho_{[\lambda t]})-[\lambda t]}{\sqrt{\lambda}}
\to_w G(k,t),
$$
where $\{G(k,t), k=1,2,\ldots,\, t\geq 0\}$ is a
mean zero Gaussian process with covariance
$$
EG(k,s)G(\ell,t)=(s\wedge t)(4(k\wedge \ell)-1_{\{k=\ell\}}-1),
$$
where weak convergence is meant on the function space $D$ that is
defined in Section {\rm 2.1} below.}

In view of Theorem J and Propositions A, B the present paper establishes
several strong approximation results in a similar vein for
random walk local times, appropriately uniformly in $k$, in both of the
cases when the time is random or deterministic.


In the next three theorems we study the asymptotic Gaussian behaviour
of the centered two-time parameter local time process
$\{\xi(k,n)-\xi(0,n),\, k=0,1,\ldots,\, n=1,2,\ldots\}$ via appropriate
strong approximations in terms of a Wiener sheet and a standard Brownian
motion.

\begin{thm}
On an appropriate probability space for a symmetric random
walk $\{S_j,\, j=0,1,\ldots\}$, we can construct a Wiener sheet
$\{W(x,y)$, $x\geq 0,y\geq 0\}$ and, independently, a standard Brownian
motion $\{W^*(y),\, y\geq 0\}$ such that, as $n\to\infty$, with
$\varepsilon>0$ we have
\begin{equation}
\xi(k,n)-\xi(0,n)=G(k,\xi(0,n))
+O(k^{5/4}n^{1/8+5\varepsilon/8})\quad {\rm a.s.}
\label{gkn}
\end{equation}
where, for a given $\varepsilon>0$, the $O(\cdot)$ term is uniform in
$k\in [1,n^{1/6-\varepsilon}]$
and
\begin{equation}
G(x,y):=W(x,y)+W(x-1,y)-W^*(y), \quad x\geq 1,\, y\geq 0.
\label{gxy}
\end{equation}
\end{thm}

The just introduced notation in (\ref{gxy}) for $G(\cdot,\cdot)$ will be
used throughout this exposition. We note that it is in fact the same
process as that of Proposition B, i.e., the two Gaussian processes agree
in distribution, but here $G(\cdot,\cdot)$ is to be constructed of
course, and so that we should have (\ref{gkn}) holding true.

It will be seen via our construction of $W(\cdot,\cdot)$ and
$W^*(\cdot)$ for establishing (\ref{gkn}) that $\xi(0,n)$ cannot be
independent of the latter Gaussian processes. This, in turn, limits
its immediate use. For the sake of making it more accessible for
applications, we also establish the next two companion conclusions to
Theorem 1.1.

For further use we introduce the notation $=_d$ for designating equality
in distribution of appropriately indicated stochastic processes.
\begin{thm}
The probability space of {\rm Theorem 1.1} can be extended to
accommodate a random walk local time $\widetilde\xi(0,n)$ such that
\begin{description}
\item{\rm (i)}\, $\{\widetilde\xi(0,n),\ n=1,2,\ldots\}=_d\{\xi(0,n),\
n=1,2,\ldots\},$
\item{\rm (ii)}\, $\widetilde\xi(0,\cdot)$ is independent of
$G(\cdot,\cdot)$
\end{description}
and, as $n\to\infty$, with $\varepsilon>0$ we have for some $\delta>0$
\begin{description}
\item{\rm (iii)}\,
$\xi(0,n)-\widetilde\xi(0,n)=O(n^{1/2-\delta}) \quad
\mathrm{a.s.},$
\item{\rm (iv)}\,
$\xi(k,n)-\xi(0,n)=G(k,\widetilde\xi(0,n))$

$+O(k^{5/4}n^{1/8+5\varepsilon/8}+kn^{1/6+\varepsilon/4}+k^{1/2}n^{1/4-\delta})
\quad \mathrm{a.s.},$
\end{description}
where, for a given $\varepsilon>0$, the latter $O(\cdot)$ term is
uniform in $k\in [1,n^{1/6-\varepsilon}]$.
\end{thm}

\begin{thm}
The probability space of {\rm Theorem 1.1} can be extended to
accommodate a standard Brownian local time process $\{\eta(0,t),\, t\geq
0\}$ such that
\begin{description}
\item{\rm (i)}\, $\eta(0,\cdot)$  is independent of
$G(\cdot,\cdot)$
\end{description}
and, as $n\to\infty$, with $\varepsilon>0$ we have for some $\delta>0$
\begin{description}
\item{\rm (ii)} $\xi(0,n)-\eta(0,n)=O(n^{1/2-\delta})\quad
\mathrm{a.s.},$
\item{\rm (iii)} $\xi(k,n)-\xi(0,n)=G(k,\eta(0,n))$

$+O(k^{5/4}n^{1/8+5\varepsilon/8}+kn^{1/6+\varepsilon/4}+
k^{1/2}n^{1/4-\delta})
\quad \mathrm{a.s.},$
\end{description}
where, for a given $\varepsilon>0$, the latter $O(\cdot)$ term is
uniform in $k\in [1,n^{1/6-\varepsilon}]$.
\end{thm}

The proofs of the above theorems will be based on the following
propositions.

\begin{propo} On an appropriate probability space for the symmetric
random walk $\{S_k,\, k=1,2,\ldots\}$ one can construct a Wiener sheet
$\{W(\cdot,\cdot)\}$ such that as $N\to\infty$, with $\varepsilon>0$ we
have
\begin{equation}
\xi(k,\rho_N^+,\uparrow)-\xi(0,\rho_N^+,\uparrow)
=W(k,2N)+O(k^{5/4}N^{1/4+\varepsilon/2})\quad {\rm a.s.},
\label{inv5}
\end{equation}
where, for a given $\varepsilon>0$, the $O$ term is uniform in
$k\in [1,N^{1/3-\varepsilon})$.
\end{propo}
\begin{propo} The probability space of {\rm Proposition 1.1} can be so
extended that as $N\to\infty$, with $\varepsilon>0$ and
$G(\cdot,\cdot)$ as in {\rm Theorem 1.1} we have
\begin{equation}
\xi(k,\rho_N^+)-\xi(0,\rho_N^+)=G(k,2N)+
O(k^{5/4}N^{1/4+\varepsilon/2})\quad {\rm a.s.},
\label{inv6}
\end{equation}
where, for a given $\varepsilon>0$, the $O$ term is uniform in
$k\in [1,N^{1/3-\varepsilon})$.
\end{propo}
\begin{propo} On the probability space of {\rm Proposition 1.1}, as
$N\to\infty$, with $\varepsilon>0$ we have
\begin{equation}
\xi(k,\rho_N)-\xi(0,\rho_N)=G(k,N)+O(k^{5/4}N^{1/4+\varepsilon/2})
\quad{\rm a.s.},
\label{inv7}
\end{equation}
where for a given $\varepsilon>0$, the $O$ term is uniform in
$k\in [1,N^{1/3-\varepsilon})$.
\end{propo}

From now on the outline of this paper is as follows. In Section 2 we
mention and prove some consequences of our just stated theorems and
propositions. In Section 3 we collect preliminary results that are
needed to prove these theorems and propositions. Theorem 1.1 and
Propositions 1.1-1.3 are proved in Section 4, while Theorems 1.2 and 1.3
in Section 5.

\section{Consequences}
\renewcommand{\thesection}{\arabic{section}} \setcounter{equation}{0}
\setcounter{thm}{0} \setcounter{lemma}{0}

Here we establish a few consequences of our theorems and propositions,
concerning weak convergence and laws of the iterated logarithm.

\subsection{Weak convergence} We start with convenient strong
approximations for the sake of concluding corresponding weak
convergence.
\begin{thm} Let $\xi(\cdot,\cdot)$, $\eta(\cdot,\cdot)$, and
$G(\cdot,\cdot)$ be as in {\rm Theorem 1.3}.
As $\lambda\to\infty$, we have
\begin{equation}
\max_{1\leq k\leq K}\, \sup_{0\leq t\leq T}
\left|\frac{\xi(k,[\lambda t])-\xi(0,[\lambda t])}{\lambda^{1/4}}
-\frac{G(k,\eta(0,\lambda t))}{\lambda^{1/4}}\right|\to 0 \quad {\rm
a.s.}
\label{weak1}
\end{equation}
and
\begin{equation}
\max_{1\leq k\leq K}\, \sup_{0\leq t\leq T}
\left|\frac{\xi(k,\rho_{[\lambda t]})-\lambda t}{\lambda^{1/2}}
-\frac{G(k,\lambda t)}{\lambda^{1/2}}\right|\to 0 \quad {\rm a.s.}
\label{weak2}
\end{equation}
for all fixed integer $K\geq 1$  and $T>0$.
\end{thm}

\noindent{\bf Proof.}
In view of Theorem 1.3 and Proposition 1.3 the respective statements of
(\ref{weak1}) and (\ref{weak2}) are seen to be true.
$\Box$

Let $N^+:=[1,2,\ldots)$, and define the space of real valued bivariate
functions
$$
f(k,t)\in D:=D(N^+\times [0,\infty))
$$
that are cadlag in $t\in [0,\infty)$. Define also
$$
\Delta=\Delta(f_1,f_2)=\max_{1\leq k\leq K}\, \sup_{0\leq t\leq T}
|f_1(k,t)-f_2(k,t)|
$$
with any fixed $(K,T)\in N^+\times [0,\infty)$, and
the measurable space $(D,{\cal D})$, where ${\cal D}$
is the $\sigma$-field generated by the $\Delta$-open balls of $D$.

On account of having for each $\lambda>0$
$$
\left\{\frac{G(k,\eta(0,\lambda t))}{\lambda^{1/4}},\,
(k,t)\in N^+\times [0,\infty)\right\}
=_d\{G(k,\eta(0,t)),\, (k,t) \in N^+\times [0,\infty)\}
$$
and
$$
\left\{\frac{G(k,\lambda t)}{\lambda^{1/2}},\,
(k,t)\in N^+\times [0,\infty)\right\}=_d
\{G(k,t),\, (k,t)\in N^+\times [0,\infty)\},
$$
Theorem 2.1 yields the following weak convergence results.
\begin{corollary}
Let $\xi(\cdot,\cdot)$, $\eta(\cdot,\cdot)$, and
$G(\cdot,\cdot)$ be as in {\rm Theorem 1.3}. As $\lambda\to\infty$,
we have
$$
h\left(\frac{\xi(k,[\lambda t])-\xi(0,[\lambda
t])}{\lambda^{1/4}}\right)
\to_d h(G(k,\eta(0,t)))
$$
and
$$
h\left(\frac{\xi(k,\rho_{[\lambda t]})-\lambda t}{\lambda^{1/2}}\right)
\to_d h(G(k,t))
$$
for all $h:D\to R$ that are $(D,{\cal D})$ measurable and
$\Delta$-continuous, or $\Delta$-continuous except at points forming a
set of measure zero on $(D,{\cal D})$ with respect to $G(\cdot,\cdot)$,
over all compact sets in $\cal D$.
\end{corollary}

\subsection{Law of the iterated logarithm}
\begin{thm} Let $K=K(t),\, t\geq 0$ be an integer valued non-decreasing
function of $t$ such that $K(t)\geq 1$ and
$$
\lim_{\alpha\to 1}\lim_{\ell\to\infty}\frac{K(\alpha^\ell)}
{K(\alpha^{\ell-1})}=1.
$$
If $K(N)\leq N^{1/3-\varepsilon}$ for some $\varepsilon>0$, then
\begin{equation}
\limsup_{N\to\infty}
\frac{\sup_{1\leq k\leq K}|\xi(k,\rho_N)-N|}
{(4K-2)^{1/2}(N\log\log N)^{1/2}}=2^{1/2}\quad
\mathrm{a.s.}
\label{strong1}
\end{equation}
If, however, $K(n)\leq n^{1/6-\varepsilon}$ for some $\varepsilon>0$,
then
\begin{equation}
\limsup_{n\to\infty}\frac{\sup_{1\leq k\leq K}|\xi(k,n)-\xi(0,n)|}
{(4K-2)^{1/2}n^{1/4}(\log\log n)^{3/4}}=\frac{2}{3}6^{1/4}\quad
\mathrm{a.s.}
\label{strong2}
\end{equation}
\end{thm}

The proof of Theorem 2.2 is based on the following result.
\begin{lemma}
For any $\alpha>1$, $K\geq 1$, $t>0$ we have the following inequalities:
\begin{equation}
P\left(\max_{1\leq k\leq K}\sup_{0\leq s\leq t}|G(k,s)|>u\right)
\label{sup1}
\end{equation}
$$
\leq C\exp\left(-\frac{u^2}{2\alpha t(4K-2)}\right),\quad u>0,
$$
\begin{equation}
P\left(\max_{1\leq k\leq K}\sup_{0\leq s\leq \eta(0,t)}|G(k,s)|
>u\right)
\label{sup2}
\end{equation}
$$
\leq
C\exp\left(-\frac{3u^{4/3}}{2^{5/3}\alpha t^{1/3}(4K-2)^{2/3}}
\right),\quad u>0
$$
with a certain positive constant $C$ depending on $\alpha$.
\end{lemma}

\noindent{\bf Proof.}
Consider the process $\{Y(s)=\max_{1\leq k\leq K}G(k,s),\, s\geq 0\}$.
$Y(s)$ is a submartingale with respect to ${\cal F}_s$, the sigma
algebra generated by $G(k,u)$, $1\leq k\leq K,\, 0\leq u\leq s$, since
if $k_0$ is defined by $G(k_0,s)=\max_{1\leq k\leq K}G(k,s)$, then obviously
$$
E(G(k_0,t)\mid {\cal F}_s)=G(k_0,s)
$$
and
$$
E(Y(t)\mid {\cal F}_s)\geq E(G(k_0,t)\mid {\cal F}_s)=G(k_0,s)=Y(s).
$$
Consequently,
$$
\left\{\sup_{0\leq s\leq t}\max_{1\leq k\leq K}G(k,s),\, t\geq
0\right\}
$$
and for $\lambda>0$
$$
\left\{\sup_{0\leq s\leq t}\max_{1\leq k\leq K}
\exp(\lambda G(k,s)), \, t\geq 0\right\}
$$
are submartingales. Using Doob inequalities (twice) we get
\begin{eqnarray}
&&P\left(\sup_{0\leq s\leq t}\max_{1\leq k\leq K}G(k,s)\geq u\right)
\nonumber\\
&&=P\left(\sup_{0\leq s\leq t}\max_{1\leq k\leq K}\exp(\lambda G(k,s))
\geq \exp(\lambda u)\right)\nonumber\\
&&\leq e^{-\lambda u}E\left(\max_{1\leq k\leq K}\exp(\lambda
G(k,t)\right)\nonumber\\
&&\leq\left(\frac{\hat\alpha}{\hat\alpha-1}\right)^{\hat\alpha}e^{-\lambda
u}E(\exp(\hat\alpha\lambda G(K,t)))
\label{doob}
\end{eqnarray}
for any $\hat\alpha>1$. $G(K,t)$ has normal distribution with mean zero
and variance $(4K-2)t$, hence
$$
E(\exp(\hat\alpha\lambda G(K,t)))=
\exp\left(\frac{\lambda^2\hat\alpha^2}{2}(4K-2)t\right)
$$
and putting
$$
\lambda=\frac{u}{\hat\alpha^2(4K-2)t},\quad \alpha=\hat\alpha^2
$$
into (\ref{doob}), we get (\ref{sup1}).

On the other hand, if $\eta(0,\cdot)$ is a Brownian local time,
independent of $G(\cdot,\cdot)$, we get from (\ref{doob})
$$
P\left(\sup_{0\leq s\leq \eta(0,t)}\max_{1\leq k\leq K}G(k,s)\geq
u\right)\leq Ce^{-\lambda u} E(\exp(\tilde\alpha\lambda
G(K,\eta(0,t)))).
$$
But
$$
\frac{G(K,\eta(0,t))}{(4K-2)^{1/2}t^{1/4}}=_d {\cal N}_1
|{\cal N}_2|^{1/2},
$$
where ${\cal N}_1$ and ${\cal N}_2$ are independent standard normal
variables. Hence
$$
P\left(\max_{1\leq k\leq K}\sup_{0\leq s\leq \eta(0,t)}G(k,s)\geq
u\right)\leq Ce^{-\lambda u}
E\left(\frac{\tilde\alpha^2\lambda^2}{2}|{\cal N}_2|(4K-2)t^{1/2}\right)
$$
$$
\leq 2Ce^{-\lambda u}
E\left(\frac{\tilde\alpha^2\lambda^2}{2}{\cal N}_2(4K-2)t^{1/2}\right)=
2Ce^{-\lambda u}
\exp\left(\frac{\tilde\alpha^4\lambda^4}{8}(4K-2)^2t\right).
$$
Putting $\lambda=(2u)^{1/3}\tilde\alpha^{-4/3}(4K-2)^{-2/3}t^{-1/3}$,
$\alpha=\tilde\alpha^{4/3}$, we get (\ref{sup2}).
$\Box$

\noindent
\textbf{Proof of Theorem 2.2.}
Let $t_\ell=\alpha^\ell,\, \, \alpha>1$. Putting
$$
u=(2\alpha^2(4K(t_\ell)-2)t_\ell\log\log t_\ell)^{1/2},\quad t=t_\ell,
K=K(t_\ell)
$$
into (\ref{sup1}), using Borel-Cantelli lemma and interpolating between
$t_{\ell-1}$ and $t_\ell$, the usual procedure gives for all large $t$
\begin{equation}
\max_{1\leq k\leq K}|G(k,t)|\leq \alpha^{3/2}
(2(4K-2)t\log\log t)^{1/2}.
\label{maxgkt}
\end{equation}
By Proposition 1.3 we also have for large $N$
$$
\max_{1\leq k\leq K}|\xi(k,\rho_N)-N|\leq
\alpha^{3/2}(2(4K-2)N\log\log N)^{1/2}+O(K^{5/4}N^{1/4+5\varepsilon/8}).
$$
Since
$$
\lim_{N\to\infty}\frac{K^{5/4}N^{1/4+\varepsilon/2}}
{(KN\log\log N)^{1/2}}=0,
$$
if $K\leq N^{1/3-\varepsilon}$, and $\alpha>1$ is arbitrary, we have an
upper bound in (\ref{strong1}).

The upper bound in (\ref{strong2}) is similar. Put
$$
u=2^{5/4}3^{-3/4}\alpha^{3/2}(4K-2)^{1/2}t^{1/4}(\log\log t)^{3/4}
$$
into (\ref{sup2}). Then, as before, we conclude that almost surely
$$
\max_{1\leq k\leq K}|G(k,\eta(0,t))|
\leq \alpha^2 2^{5/4}3^{-3/4}(4K-2)^{1/2}t^{1/4}(\log\log t)^{3/4}
$$
for $t$ large enough. Since $\alpha>1$ is arbitrary, using Theorem 1.3,
we get an upper bound in (\ref{strong2}).

To prove the lower bound in (\ref{strong1}), for $0<\delta<1$ define
the events
$$
A_\ell=\{G(K_\ell,t_\ell)-G(K_\ell,t_{\ell-1})\geq
(1-\delta)(2(4K_\ell-2)t_\ell\log\log t_\ell)^{1/2}\},
$$
$\ell=1,2,\ldots$, where $t_\ell=\delta^{-\ell}$ and $K_\ell=K(t_\ell)$.
Since $G(K_\ell,t_\ell)-G(K_\ell,t_{\ell-1})$ has normal distribution
with mean zero and variance $(4K_\ell-2)(t_\ell-t_{\ell-1})$, an easy
calculation shows
$$
P(A_\ell)\geq \frac{C}{(\log
t_\ell)^{1-\delta}}=\frac{C}{\ell^{1-\delta}}.
$$
Since $A_\ell$ are independent, Borel-Cantelli lemma implies
$P(A_\ell\, \mathrm{i.o.})=1$. But
$$
G(K_\ell,t_{\ell-1})\leq (1+\delta)(4K_\ell-2)^{1/2}
(2t_\ell\log\log t_\ell)^{1/2}\delta^{1/2}
$$
for all large $\ell$, we have also
$$
G(K_\ell,t_\ell)\geq ((1-\delta)-(1+\delta)\delta^{1/2})
(2(4K_\ell-2)t_\ell\log\log t_\ell)^{1/2}
$$
infinitely often with probability 1. Since $\delta>0$ is arbitrary, we
conclude
$$
\limsup_{t\to\infty}\frac{G(K,t)}{(2(4K-2)t\log\log t)^{1/2}}\geq 1.
$$
Using Proposition 1.3, this also gives a lower bound in (\ref{strong1}).

To show the lower bound in (\ref{strong2}), we follow Burdzy \cite{bu}
with some modifications. Define $t_\ell=\exp(\ell\log\ell)$, and the
events
$$
A_\ell^{(1)}=\{(1-\delta)a_\ell\leq \eta(t_\ell)\leq
2(1-\delta)a_\ell\}
$$
and
$$
A_\ell^{(2)}=\{\inf_{s\in I_\ell}G(K_\ell,s)-G(K_\ell,\gamma a_\ell)
\geq (1-2\beta)(4K_\ell-2)^{1/2}u_\ell\},
$$
where $K_\ell=K(t_\ell)$, $\eta(t_\ell)=\eta(0,t_\ell)$,
$$
a_\ell=\left(\frac{2}{3}t_\ell\log\log t_\ell\right)^{1/2},
$$
$$
u_\ell=\frac{2}{3^{1/2}}((1-2\delta)a_\ell\log\log a_\ell)^{1/2},
$$
$$
I_\ell=[(1-2\delta)a_\ell, 3(1-2\delta)a_\ell],
$$
and $\beta,\, \delta,\, \gamma$ are certain small constants to be
choosen later on. Obviously, the events $\{A_\ell^{(1)},\,
\ell=1,2,\ldots\}$ and
$\{A_\ell^{(2)},\, \ell=1,2,\ldots\}$ are independent. Let
$$
A_\ell=A_\ell^{(1)} A_\ell^{(2)}.
$$
We show that for certain values of the above constants, $P(A_\ell\,
\mathrm{i.o.})=1$.

Since $t_\ell^{-1/2}\eta(t_\ell)$
is distributed as the absolute value of a standard normal variable, an
easy calculation shows
$$
P(A_\ell^{(1)})\geq \frac{C}{(\log t_\ell)^{(1-\delta)/3}}
$$
with some $C>0$.

Converting the inequality of Lemma 2 in Burdzy \cite{bu} from small time
to large time, the following inequality can be concluded for large
enough $u$:
$$
P(\inf_{s\in [u,3u]}W(s)-W(\gamma u)\geq (1-2\beta)(2/3^{1/2})
(u\log\log u)^{1/2})
$$
$$
\geq (\log u)^{(-2/3)(1-\beta)/(1-\gamma)},
$$
where $W(\cdot)$ is a standard Wiener process. From this we get
$$
P(A_\ell^{(2)})\geq \frac{C}{(\log u_\ell)^{2(1-\beta)/(3(1-\gamma))}}
$$
with $C>0$.

We can choose the constants $\beta,\gamma,\delta$ appropriately to have
$\sum_{\ell}P(A_\ell)=\sum_\ell P(A_\ell^{(1)})P(A_\ell^{(2)})=\infty$.
The events $A_\ell$ however are not independent. Next we show
$P(A_jA_\ell)\leq CP(A_j)P(A_\ell)$ with some constant $C$.
It can be seen that for large $\ell$ we have
$3(1-2\delta)a_\ell\leq \gamma a_{\ell+1}$, therefore $A_\ell^{(2)}$ are
independent events for $\ell\geq \ell_0$ with a certain $\ell_0$. We
have
$$
P(A_jA_\ell)=P(A_j^{(1)}A_\ell^{(1)})P(A_j^{(2)})P(A_\ell^{(2}).
$$
It suffices to show that $P(A_j^{(1)}A_\ell^{(1)})\leq
CP(A_j^{(1)})P(A_\ell^{(1)})$.
For this purpose it is more convenient to work with $M(t)$, the supremum
of the Wiener process, since according to L\'evy's theorem,
the process $\{\eta(0,t),\, t\geq 0\}$ is identical in distribution with
$\{M(t),\, t\geq 0\}$. So let $\{\widehat W(t),\, t\geq 0\}$ be a
standard Wiener process and $M(t)=\sup_{0\leq s\leq t}\widehat W(s)$.
Denote by $g_t(y)$ the density of $M(t)$ and by $g_{t_1,t_2}(y_1,y_2)$
the joint density of $M(t_1), M(t_2)$ . It is well known that
$$
g_t(y)=\frac{2}{\sqrt{2\pi t}}\exp\left(-\frac{y^2}{2t}\right).
$$
Then with $h(t_1,z)$, the joint density of $M(t_1)$ and
$\widehat W(t_1)$, we can write for $t_1<t_2$,
$$
g_{t_1,t_2}(y_1,y_2)=\int_{-\infty}^{y_1} h(t_1,z)g_{t_2-t_1}(y_2-z)\,
dz.
$$
It can be seen that for $z\leq y_1$
$$
g_{t_2-t_1}(y_2-z)\leq \sqrt{\frac{t_2}{t_2-t_1}}g_{t_2}(y_2)
\exp\left(\frac{y_1y_2}{t_2-t_1}\right).
$$
Hence
$$
g_{t_1,t_2}(y_1,y_2)\leq
g_{t_1}(y_1)g_{t_2}(y_2)\sqrt{\frac{t_2}{t_2-t_1}}
\exp\left(\frac{y_1y_2}{t_2-t_1}\right).
$$
Returning to the probability of the events $A^{(1)}$, we have for
$j<\ell$
$$
P(A_j^{(1)} A_\ell^{(1)})\leq\sqrt{\frac{t_\ell}{t_\ell-t_j}}
\exp\left(\frac{4(1-\delta)^2a_ja_\ell}{t_\ell-t_j}\right)
P(A_j^{(1)})P(A_\ell^{(1)})
$$
$$
\leq CP(A_j^{(1)})P(A_\ell^{(1)}),
$$
where $C>1$ can be choosen arbitrarily close to $1$ by choosing $\ell-j$
sufficiently large. Hence for any $\varepsilon>0$ there exists $m_0$
such that
$$
P(A_jA_\ell)\leq (1+\varepsilon)P(A_j)P(A_\ell)
$$
if $\ell-j\geq m_0$. It follows that
$$
\sum_{\ell=1}^n\sum_{j=1}^\ell P(A_j A_\ell)
\leq (1+\varepsilon)\sum_{\ell=1}^n\sum_{j=1}^{\ell-m_0} P(A_j)P(A_\ell)
+m_0\sum_{\ell=1}^n P(A_\ell).
$$
By Borel-Cantelli lemma (cf. \cite{sp}, p. 317) $P(A_\ell\,
\mathrm{i.o.})\geq 1/(1+\varepsilon)$. Since $\varepsilon>0$ is
arbitrary, we also have $P(A_\ell\, \mathrm{i.o.})= 1$.
$A_\ell^{(1)}$ implies
$$
\eta(t_\ell)\in [(1-\delta)a_\ell,2(1-\delta) a_\ell]
\subset I_\ell,
$$
consequently, if both $A_\ell^{(1)}$ and $A_\ell^{(2)}$ occur, then
$$
G(K_\ell,\eta(t_\ell))
$$
$$
\geq (1-2\beta)(4K_\ell-2)^{1/2}2
((1-2\delta)a_\ell\log\log a_\ell)^{1/2}/3^{1/2}+G(K_\ell,\gamma
a_\ell).
$$
It follows from (\ref{maxgkt}) that
$$
G(K_\ell,\gamma a_\ell)
\geq -(4K_\ell-2)^{1/2}(\gamma a_\ell\log\log
a_\ell)^{1/2}
$$
for all large $\ell$ with probability 1, i.e.
$$
\limsup_{\ell\to\infty}\frac{G(K_\ell,\eta(t_\ell))}
{(4K_\ell-2)^{1/2}(a_\ell\log\log a_\ell)^{1/2}}
\geq (1-2\beta)(1-2\delta)^{1/2}2/3^{1/2}
-\gamma^{1/2}.
$$
But
$$
\lim_{\ell\to\infty}\frac{a_\ell\log\log a_\ell}
{t_\ell^{1/2}(\log\log t_\ell)^{3/2}}=(2/3)^{1/2},
$$
implying
$$
\limsup_{t\to\infty}\frac{\max_{1\leq k\leq K}G(k,\eta(0,t))}
{(4K-2)^{1/2}t^{1/4}(\log\log t)^{3/4}}
$$
$$
\geq 2^{5/4}3^{-3/4}
(1-2\beta)(1-2\delta)^{1/2}-2^{1/4}3^{-1/4}\gamma^{1/2}
\quad\mathrm{a.s.}
$$
Since it is possible to choose $\beta, \delta, \gamma$ arbitrarily
small, combining this with Theorem 1.3, gives a lower bound in
(\ref{strong2}).
$\Box$

\section{Preliminaries}
\renewcommand{\thesection}{\arabic{section}} \setcounter{equation}{0}
\setcounter{thm}{0} \setcounter{lemma}{0}

In this Section we collect the results needed to prove our theorems and
propositions. The proofs will use the branching property (Ray-Knight
description) of the random walk local time. For more details in this
respect we refer to Knight \cite{kn}, Dwass \cite{dw}, Rogers \cite{ro}
and T\'oth \cite{to}.

Introduce the following notations for $k=1,2,\ldots,\, i=1,2,\ldots$.
\begin{equation}
\tau_i^{(k)}:=\min\{j>\tau_{i-1}^{(k)}:\, S_{j-1}=k, S_j=k-1\},
\label{tauik}
\end{equation}
with $\tau_0^{(k)}:=0$,
\begin{equation}
T_i^{(k)}:=\xi(k,\tau_i^{(k)})-\xi(k,\tau_{i-1}^{(k)}).
\label{tik}
\end{equation}
With probability 1, there is such a double infinite sequence of
$\tau_i^{(k)}$ and hence also of $T_i^{(k)},\, i=1,2,\ldots,\,
k=1,2,\ldots$
\begin{lemma}
The random variables $\{T_i^{(k)},\, k=1,2,\ldots, i=1,2,\ldots\}$ are
completely independent and distributed as
\begin{equation}
P(T_i^{(k)}=j)=\frac{1}{2^j},\quad j=1,2,\ldots,
\label{tik2}
\end{equation}
\begin{equation}
E(T_i^{(k)})=2,\quad Var(T_i^{(k)})=2.
\end{equation}
\end{lemma}

\noindent{\bf Proof.}
Obvious.

Introduce
\begin{equation}
U^{(k)}(j):=T_1^{(k)}+\ldots +T_j^{(k)}-2j, \quad k=1,2,\ldots,\, \,
j=1,2,\ldots\, .
\label{ukj}
\end{equation}

For the following inequality we refer to T\'oth \cite{to}.
\begin{lemma}
$$
P\left(\max_{1\leq i\leq n}|U^{(k)}(i)|>z\right)
\leq 2\exp\left(-\frac{z^2}{8n}\right),\quad 0<z<a_0n
$$
for some $a_0>0$.
\end{lemma}

We need Hoeffding's inequality \cite{ho} for binomial distribution
(cf. also Sho\-rack and Wellner \cite{sw}, pp. 440).
\begin{lemma}
Let $\nu_N$ have binomial distribution with parameters $(N,1/2)$. Then
$$
P(|2\nu_N-N|\geq u)\leq 2\exp\left(-\frac{u^2}{2N}\right),\quad
0<u.
$$
\end{lemma}

To establish our results, we make use of one of the celebrated KMT
strong invariance principles (cf. Koml\'os et al. \cite{kmt}).

\begin{lemma} Let $\{Y_i\}_{i=1}^\infty$ be i.i.d. random variables with
expectation zero, variance $\sigma^2$ and having moment generating
function in a neighbourhood of zero. On an appropriate probability space
one can construct $\{Y_i\}_{i=1}^{\infty}$  and a Wiener process $\{W(t),\,
t\geq 0\}$ such that for all $x>0$ and $n=1,2,\ldots$
$$
P\left(\max_{1\leq i\leq n}\left|\sum_{j=1}^i
Y_j-W(i\sigma^2)\right|>C_1\log n+x\right)\leq C_2e^{-C_3 x},
$$
where $C_1,C_2,C_3$ are positive constants, and $C_3$ can be choosen
arbitrarily large by choosing $C_1$ sufficiently large.
\end{lemma}

There are several papers on strong invariance principles for local
times, initiated by R\'ev\'esz \cite{re1}, and further developed by
Borodin \cite{bo2}, \cite{bo}, Bass and Khoshnevisan \cite{bk}, and
others, as in the references of these papers. The best rate via
R\'ev\'esz's Skorokhod type construction was given by Cs\"org\H o and
Horv\'ath \cite{csh}.
\begin{lemma}
On a rich enough probability space one can define a simple symmetric
random walk with local time $\xi(\cdot,\cdot)$ and a standard Brownian
local time $\eta(\cdot,\cdot)$ such that as $n\to\infty$
\begin{equation}
\sup_{k\in Z}|\xi(k,n)-\eta(k,n)|=
O(n^{1/4}(\log n)^{1/2}(\log\log n)^{1/4})\quad {\rm a.s.}
\label{locinv}
\end{equation}
\end{lemma}

We note in passing that having (\ref{locinv}) with
$O(n^{1/4}(\log\log n)^{3/4})$ is best possible for any construction
(cf. Cs\"org\H o and Horv\'ath \cite{csh}), i.e., only the $(\log
n)^{1/2}$ term of (\ref{locinv}) could be changed, and only to
$(\log\log n)^{1/2}$, by any other construction. It remains an open
problem to find such a construction that would achieve this best
possible minimal gain.

\begin{lemma}
Let $\{W_i(\cdot),\, i=1,\ldots,k\}$ be independent Wiener
processes and $t>0$. The following inequality holds.
$$
P\left(\sup_{0\leq t_i\leq t,\, i=1,\ldots,k}\left|\sum_{i=1}^k
W_i(t_i)\right|\geq z\right)\leq 2k e^{-z^2/(2k^2t)},\quad 0<z.
$$
\end{lemma}

\noindent{\bf Proof.}
Since
$$
\sup_{0\leq t_i\leq t,\, i=1,\ldots,k}\left|\sum_{i=1}^k W_i(t_i)\right|
\leq k\max_{1\leq i\leq k}\sup_{0\leq t_i\leq t}\left|W_i(t_i)\right|,
$$
we have
$$
P\left(\sup_{0\leq t_i\leq t,\, i=1,\ldots,k}\left|\sum_{i=1}^k
W_i(t_i)\right|\geq z\right)
\leq P\left(\max_{1\leq i\leq k}\sup_{0\leq t_i\leq
t}\left|W_i(t_i)\right|\geq z/k\right)
$$
$$
\leq kP\left(\sup_{0\leq s\leq t}\left|W(s)\right|\geq z/k\right)\leq
2k e^{-z^2/(2k^2t)}.
$$
$\Box$

\begin{lemma} The following identities hold.
\begin{equation}
\xi(k,\rho_N^+)=U^{(k)}(\xi(k-1,\rho_N^+,\uparrow))+
2\xi(k-1,\rho_N^+,\uparrow),
\end{equation}
\begin{equation}
\xi(k,\rho_N^+,\uparrow)=\sum_{i=1}^{\xi(k-1,\rho_N^+,\uparrow)}
(T_i^{(k)}-1),
\label{branching}
\end{equation}
\begin{equation}
\xi(k,\rho_N^+,\downarrow)=\xi(k-1,\rho_N^+,\uparrow).
\label{xidown}
\end{equation}
\end{lemma}

\noindent{\bf Proof.}
Obvious.

Equation (\ref{branching}) amounts to saying that
$\xi(k,\rho_N^+,\uparrow),\, k=0,1,\ldots$ is a critical branching
process with geometric offspring distribution.

\begin{lemma} For $K\geq 1$
\begin{equation}
P(\max_{1\leq k\leq K}\xi(k,\rho_N^+)\geq 5N)\leq
K\exp\left(-\frac{N}{4K}\right)
\end{equation}
\end{lemma}

\noindent{\bf Proof.}
For the distribution of $\xi(k,\rho_1^+)$ we have (cf. R\'ev\'esz
\cite{re2})
\begin{equation}
P(\xi(k,\rho_1^+)=m)=\left\{
\begin{array}{ll}
& 1-\frac{1}{k}\quad\quad\quad\quad\quad\quad\quad \mathrm{if}\quad
m=0,\\
&\frac{1}{2k^2}\left(1-\frac{1}{2k}\right)^{m-1}\quad\quad
\mathrm{if}\quad
m=1,2,\ldots
\end{array}
\right.
\end{equation}
Hence for the moment generating function we have
$$
g(k,t)=E\left(e^{t\xi(k,\rho_1^+)}\right)=
\frac{1-(2k-2)(1-e^{-t})}{1-2k(1-e^{-t})}\leq 1+\frac{2t}{1-2kt}.
$$
Selecting $t=1/(4K)$, we arrive at
$$
g(k,t)\leq 1+\frac{1}{2K-k}\leq 1+\frac1K\leq e^{1/K}.
$$
Since
$$
E\left(e^{t\xi(k,\rho_N^+)}\right)=(g(k,t))^N,
$$
we have by Markov's inequality
$$
P(\xi(k,\rho_N^+)\geq 5N)\leq\left(g(k,t)e^{-5t}\right)^N
\leq e^{\left(\frac1K-\frac5{4K}\right)N}= e^{-N/(4K)}.
$$
$\Box$

We need inequalities for increments of the Wiener process (Cs\"org\H o
and R\'ev\'esz \cite{csr}), Brownian local time (Cs\'aki et al.
\cite{cscsfr}), and random walk local time (Cs\'aki and F\"oldes
\cite{csf83}).
\begin{lemma}
With any constant $C_2<1/2$ and some $C_1>0$ we have
$$
P\left(\sup_{0\leq s\leq T-h}\sup_{0\leq t\leq h}|W(s+t)-W(s)|
\geq v\sqrt{h}\right)
\leq \frac{C_1T}{h}e^{-C_2v^2},
$$
$$
P\left(\sup_{0\leq s\leq t-h}(\eta(0,h+s)-\eta(0,s))\geq
x\sqrt{h}\right)
\leq C_1\left(\frac{t}{h}\right)^{1/2}e^{-C_2x^2},
$$
and
$$
P\left(\max_{0\leq j\leq t-a}(\xi(0,a+j)-\xi(0,j))\geq x\sqrt{a}\right)
\leq C_1\left(\frac{t}{a}\right)^{1/2}e^{-C_2x^2}.
$$
\end{lemma}
Note that we may have the same constants $C_1,C_2$ in the above
inequalities. In fact, in our proofs the values of these constants
are not important, and it is indifferent whether they are the same
or not. We continue using these notations for constants of no interest
that may differ from line to line.

\begin{lemma} For $1\leq u$ we have
$$
P(\rho_N\geq uN^2)\leq \frac{1}{\sqrt{u}}
$$
and
$$E(\rho_1 I\{\rho_1\leq u\})\leq 3 \sqrt{u}.$$
\end{lemma}

\noindent{\bf Proof.}
For the distribution of $\rho_N$ we have (cf. R\'ev\'esz \cite{re2}, pp.
98)
$$
P(\rho_N>2n)=\frac{1}{2^{2n}}\sum_{j=0}^{N-1}2^j{2n-j\choose n},
\quad n=1,2,\ldots.
$$
An elementary calculation shows that the largest term in the sum above
is
for $j=0$, hence
$$
P(\rho_N>2n)\leq \frac{N}{2^{2n}}{2n\choose n}.
$$
Moreover, it can be easily seen that
$$
\frac{(2n+2)^{1/2}}{2^{2n}}{2n\choose n}
$$
is decreasing in $n=1,2,\ldots$, hence it is less than 1 for all $n$,
and thus implying
$$
P(\rho_N>2n)\leq\frac{N}{(2n+2)^{1/2}}.
$$
For a given $u\geq 1$ choose $n$ so that $2n<uN^2\leq 2n+2$. Then
$$
P(\rho_N\geq uN^2)\leq P(\rho_N>2n)\leq\frac{N}{(2n+2)^{1/2}}
\leq \frac1{\sqrt{u}}.
$$
Moreover,
$$
E(\rho_1 I\{\rho_1\leq u\})=\sum_{1\leq j\leq u}jP(\rho_1=j)
\leq \sum_{0\leq j\leq u}P(\rho_1\geq u)
$$
$$
\leq 1+\sum_{1\leq j\leq
u}\frac1{\sqrt{j}}\leq 1+\int_0^u\frac{dx}{\sqrt{x}}=1+2\sqrt{u}
\leq 3 \sqrt{u}.$$
$\Box$

\begin{lemma} Define $\tau_0:=0$,
$$
\tau_n:=\inf\{t:\, t>\tau_{n-1},\, |W(t)-W(\tau_{n-1})|=1\},
\quad n=1,2,\ldots
$$
Then $\tau_n$ is a sum of $n$ i.i.d. random variables, $E(\tau_1)=1$ and
\begin{equation}
E(e^{\theta \tau_1})=\frac{1}{\cosh(\sqrt{2\theta})}.
\label{tau1}
\end{equation}
Moreover,
\begin{equation}
P(|\tau_n-n)|\geq u\sqrt{n})\leq 2e^{-3u^2 /8}, \qquad 0<u<2\sqrt{n}/3.
\label{tau2}
\end{equation}
\end{lemma}

\noindent{\bf Proof.}
For (\ref{tau1}) see, e.g., Borodin and Salminen \cite{bs}.
To show (\ref{tau2}), we use exponential Markov's inequality:
$$
P(|\tau_n-n|\geq u\sqrt{n})
\leq e^{-u\theta\sqrt{n}}\left((g(\theta))^n+(g(-\theta))^n\right),
$$
for $0<\theta\leq 1/2$, where
$$
g(\theta):=E(e^{\theta(\tau_1-1)})=\frac1{e^\theta\cosh(\sqrt{2\theta})}.
$$
By the series expansion of $\log\cos x$ (cf. Abramowitz and Stegun
\cite{as}, pp. 75) and putting $\cosh x=\cos (ix)$, we get
$$
\log \cosh x=\sum_{k=1}^\infty \frac{2^{2k-1}(2^{2k}-1)B_{2k}}
{k(2k)!}x^{2k},\quad |x|\leq \frac{\pi}{2},
$$
where $B_i$ are Bernoulli numbers, and using that $B_2=1/6$,
$B_4=-1/30$ and the inequality (cf. \cite{as}, pp. 805)
$$
|B_{2n}|\leq\frac{2(2n)!}{(2\pi)^{2n}(1-2^{1-2n})}
$$
for $n>2$, one can easily see that
$$
\log g(\theta)\leq \frac{\theta^2}3 (1+\theta+\theta^2+\ldots)
=\frac{\theta^2}{3(1-\theta)}\leq 2\theta^2/3
$$
if $0\leq\theta\leq 1/2$. Similarly,
$$
\log g(-\theta)\leq 2\theta^2/3, \quad 0\leq\theta\leq 1/2,
$$
hence putting $\theta=3u/(4\sqrt{n})$, we get (\ref{tau2}).
$\Box$

\begin{lemma}
Let $Y_i,\, i=1,2,\ldots$ be i.i.d. random variables having
exponential distribution with parameter 1. Then
\begin{equation}
P\left(\max_{1\leq j\leq n}\left|\sum_{i=0}^j (Y_i-1)\right|\geq
u\sqrt{n}\right)\leq 2e^{-u^2/8},\qquad 0<u<2\sqrt{n}.
\label{exp}
\end{equation}
Moreover, with any $C>0$,
\begin{equation}
P\left(\max_{1\leq i\leq n} Y_i\geq C\log n\right)\leq n^{1-C}.
\label{exp2}
\end{equation}
\end{lemma}

\noindent{\bf Proof.}
By exponential Kolmogorov's inequality (see T\'oth \cite{to}) we have
for $0<\theta\leq 1/2$
$$
P\left(\max_{1\leq j\leq n}\left|\sum_{i=0}^j (Y_i-1)\right|\geq
u\sqrt{n}\right)
$$
$$
\leq e^{-\theta u\sqrt{n}}\left((f(\theta))^n+(f(-\theta))^n\right),
$$
where
$$
f(\theta)=E\left(e^{\theta(Y_1-1)}\right)=\frac{1}{e^\theta(1-\theta)}
\leq e^{2\theta^2}
$$
and, similarly,
$$
f(-\theta)=\frac{e^\theta}{1+\theta}\leq e^{2\theta^2}.
$$
Now (\ref{exp}) can be obtained by putting $\theta=u/(4\sqrt{n})$, and
(\ref{exp2}) is easily seen as follows.
$$
P\left(\max_{1\leq i\leq N}Y_i\geq C\log N\right)\leq NP(Y_1\geq C\log
N)=N^{1-C}.
$$
$\Box$

Finally, we quote the following lemma from Berkes and Philipp \cite{BP}.

\begin{lemma}
Let $B_i,\, i=1,2,3$ be separable Banach
spaces. Let $F$ be a distribution on $B_1\times B_2$ and let $G$
be a distribution on $B_2\times B_3$ such that the second marginal
of $F$ equals the first marginal of $G$. Then there exists a
probability space and three random variables $Z_i,\, i=1,2,3$,
defined on it such that the joint distribution of $Z_1$ and $Z_2$
is $F$ and the joint distribution of $Z_2$ and $Z_3$ is $G$.
\end{lemma}

\section{Proof of Theorem 1.1}
\renewcommand{\thesection}{\arabic{section}} \setcounter{equation}{0}
\setcounter{thm}{0} \setcounter{lemma}{0}

\subsection{Proof of Proposition 1.1}

First we prove the next lemma, which is a consequence of Lemma 3.4.
\begin{lemma}
On an appropriate probability space one can construct independent random
variables $\{T_i^{(k)}\}_{i,k=1}^{\infty}$ with distribution {\rm
(\ref{tik2})} and a sequence of independent Wiener processes
$\{W_k(t),\, t\geq 0\}_{k=1}^\infty$ such that, as $N\to\infty$, we have
\begin{equation}
\max_{1\leq k\leq N}\max_{1\leq j\leq N}|U^{(k)}(j)-W_k(2j)|
=O(\log N)\quad {\rm a.s.},
\label{ukwk}
\end{equation}
where $U^{(k)}(j)$ are defined by {\rm (\ref{ukj})}.
\end{lemma}

\noindent{\bf Proof.}
By Lemma 3.4 for each fixed $k=1,2,\ldots$ on a probability space
one can construct $T^{(k)}_j$ and $W_k$ satisfying
\begin{equation}
P\left(\max_{1\leq j\leq N}|U^{(k)}(j)-W_k(2j)|\geq (C_1+1)\log
N\right) \leq C_2e^{-C_3\log N}. \label{kmt}
\end{equation}
Note that the constants $C_1,C_2,C_3$ depend only on the distribution of
$T^{(k)}_j$, hence they do not depend on $k$.
Now consider the product space so that we have (\ref{kmt}) for all
$k=1,2,\ldots$ on it. Then
$$
P\left(\max_{1\leq k\leq N}\max_{1\leq j\leq N}|U^{(k)}(j)-W_k(2j)|
\geq (C_1+1)\log N\right)
$$
$$
\leq NC_2e^{-C_3\log N}=C_2e^{-(C_3-1)\log N}=\frac{C_2}{N^{C_3-1}}.
$$
Choosing $C_3>2$, (\ref{ukwk}) follows by Borel-Cantelli lemma.
$\Box$

Now on the probability space of Lemma 4.1 a Wiener sheet
$W(\cdot,\cdot)$ is constructed from the
independent Wiener processes $W_k,\, k=1,2,\ldots$ as above in such a
way that for integer $k$ we have (cf. Section 1.11 of \cite{csr})
\begin{equation}
W(k,y)=\sum_{i=1}^k W_i(y).
\label{wky}
\end{equation}

By Lemma 3.13 this can be extended to a Wiener sheet
$\{W(x,y),\, x,y\geq 0\}$ on the probability space of Lemma 4.1,
so that on the same probability space we have
a simple symmetric random walk $\{S_i\}_{i=0}^\infty$ as
defined in the Introduction, satisfying (\ref{tauik}) and (\ref{tik}).

To show Proposition 1.1, we start from the identity
$$
\xi(k,\rho_N^+,\uparrow)=\xi(k,\rho_N^+)-\xi(k,\rho_N^+,\downarrow)
=U^{(k)}(\xi(k-1,\rho_N^+,\uparrow))+\xi(k-1,\rho_N^+,\uparrow).
$$
Repeating this procedure several times, we arrive at
$$
\xi(k,\rho_N^+,\uparrow)=\sum_{i=1}^k
U^{(i)}(\xi(i-1,\rho_N^+,\uparrow)) +N.
$$
For brevity, from here on in this proof we use the notation
$$
\xi_i= \xi(i,\rho_N^+,\uparrow).
$$
Continuing accordingly, using Lemma 4.1 and the fact that as
$N\to\infty$
$$
\max_{1\leq i\leq N^{1-\varepsilon}}\log \xi_i=O(\log N)\quad \qquad
{\rm a.s.},
$$
which follows from Lemma 3.8, we get for $k=1,2,\ldots$
$$
\xi_k=\sum_{i=1}^k W_i(2\xi_{i-1})+N +O(k\log N)
$$
$$
=W(k,2N)+N+O(k\log N)+\sum_{i=1}^k(W_i(2\xi_{i-1})-W_i(2N))
\quad{\rm a.s.}
$$
Now we are to estimate the last term in our next lemma.
\begin{lemma} As $N\to\infty$,
$$
\sum_{i=1}^k(W_i(2\xi_{i-1})-W_i(2N))=O(k^{5/4}N^{1/4+\varepsilon/2})
\quad\mathrm{a.s.},
$$
where the $O$ term is uniform in $k\in [1,N^{1/3-\varepsilon}]$.
\end{lemma}

\noindent{\bf Proof.}
Observe that
$$
\sum_{i=1}^k(W_i(2\xi_{i-1})-W_i(2N))=
\sum_{i=1}^k\widetilde W_i(2|\xi_{i-1}-N|),
$$
where $\widetilde W_i(\cdot),\, i=1,2,\ldots$ are independent Wiener
processes.

Let $K=[N^{1/3-\varepsilon}]$, $w_k=k^{1/2}N^{1/2+\varepsilon/2}$,
$z_k=k^{5/4}N^{1/4+\varepsilon/2}$.
Then
$$
P\left(\bigcup_{k=1}^K
\left\{\left|\sum_{i=1}^k(W_i(2\xi_{i-1})-W_i(2N))\right|\geq
z_k\right\}\right)
$$
$$
=P\left(\bigcup_{k=1}^K \left\{\left|\sum_{i=1}^k\widetilde
W_i(2|\xi_{i-1}-N|)\right|\geq z_k \right\}\right)
$$
$$
\leq \sum_{k=1}^KP\left(\left|\sum_{i=1}^k\widetilde
W_i(2|\xi_{i-1}-N|)\right|\geq z_k\right)
$$
$$
\leq\sum_{k=1}^K\left(P\left(\max_{1\leq i\leq k}|\xi_i-N|\geq
w_k\right) +P\left(\sup_{0\leq t_i\leq 2w_k,\,
i=1,\ldots,k}\left|\sum_{i=1}^k \widetilde W_i(t_i)\right|\geq
z_k\right)\right).
$$
It follows from (\ref{branching}) by telescoping that
$$
\xi_i-N=\sum_{j=1}^{\xi_{i-1}+\xi_{i-2}+\ldots \xi_1+N} (T_j-2),
\quad i=1,2,\ldots
$$
where $T_j$ are i.i.d. random variables distributed as $T^{(k)}_i$.
From Lemma 3.2 and Lemma 3.8 we obtain
\begin{eqnarray*}
&&P\left(\max_{1\leq i\leq k}|\xi_i-N|\geq w_k\right)\\
&&\leq P\left(\max_{1\leq i\leq k}\xi_i\geq 5N\right)
+kP\left(\max_{1\leq n\leq 5Nk} \left|\sum_{j=1}^n(T_j-2)\right|\geq
w_k\right)\\
&&\leq ke^{-N/(4k)}+2ke^{-w_k^2/(40kN)}.
\end{eqnarray*}
From this, together with Lemma 3.6, we finally get
$$
P\left(\bigcup_{k=1}^K
 \left\{\left|\sum_{i=1}^k(W_i(2\xi_{i-1})-W_i(2N))\right|\geq z_k \right\}\right)
$$
$$
\leq\sum_{k=1}^K\left(ke^{-N/(4k)}+2ke^{-w_k^2/(40kN)}+
2ke^{-z_k^2/(4k^2w_k)} \right)
$$
$$
\leq N^{2/3}e^{-N^{2/3+\varepsilon}/4}+
2N^{2/3}e^{-N^{\varepsilon}/40}+ 2N^{2/3}e^{-N^{\varepsilon/2}/4}.
$$
This is summable in $N$, so the lemma follows by Borel-Cantelli lemma.
$\Box$

Since $\xi(0,\rho_N^+,\uparrow)=N$, this also proves Proposition 1.1.
$\Box$

\subsection{Proof of Proposition 1.2}

According to (\ref{xidown}) and Proposition 1.1, as $N\to\infty$,
\begin{eqnarray}
&&\xi(k,\rho_N^+)=\xi(k,\rho_N^+,\uparrow)+\xi(k-1,\rho_N^+,\uparrow)
\nonumber\\
&&=2N+W(k,2N)+W(k-1,2N)+O(k^{5/4}N^{1/4+\varepsilon/2})\quad {\rm a.s.}
\label{xi+}
\end{eqnarray}

On the other hand,
\begin{equation}
\xi(0,\rho_N^+)=\xi(0,\rho_N^+,\uparrow)+\xi(0,\rho_N^+,\downarrow)
=N+\xi(0,\rho_N^+,\downarrow).
\label{xi+2}
\end{equation}
But
\begin{equation}
\xi(0,\rho_N^+,\downarrow)=T_1^*+\ldots +T_N^*,
\label{xi+down}
\end{equation}
where $T_i^*$ represents the number of downward excursions away from $0$
between the $i$th and $(i+1)$st upward excursions away from $0$. Hence
$T_i^*$ are i.i.d. random variables with geometric distribution
$$
P(T_i^*=j)=\frac{1}{2^{j+1}},\quad j=0,1,2,\ldots
$$
and also independent of $\{T^{(k)}_i,\, i,k=1,2,\ldots\}$. Hence from
KMT Lemma 3.4 and by Lemma 3.13, on the probability space of
Proposition 1.1 one can construct a Wiener process $W^*(\cdot)$,
independent of $W(\cdot,\cdot)$ such that, as $N\to\infty$,
$$
T_1^*+\ldots T_N^*=N+W^*(2N)+O(\log N)\quad {\rm a.s.}
$$
This together with (\ref{xi+}), (\ref{xi+2}) and (\ref{xi+down}) proves
Proposition 1.2.
$\Box$

\subsection{Proof of Proposition 1.3}

Consider $N$ excursions away from $0$, out
of which $\nu_N$ are upward excursions, and $N-\nu_N$ are downward
excursions. According to Proposition 1.2, as $N\to\infty$,
$$
\xi(k,\rho_{\nu_N}^+)-\xi(0,\rho_{\nu_N}^+)=
G(k,2\nu_N)+O(k^{5/4}(\nu_N)^{1/4+\varepsilon/2})
$$
$$
=G(k,2\nu_N)+O(k^{5/4}N^{1/4+\varepsilon/2})
\quad {\rm a.s.}
$$

Since $\xi(k,\rho_N)=\xi(k,\rho_{\nu_N}^+)$ for $k>0$, it is enough
to verify the next Lemma.
\begin{lemma} As $N\to\infty$ we have
\begin{equation}
W(k,2\nu_N)-W(k,N)=O(k^{1/2}N^{1/4+\varepsilon/2})\quad {\rm a.s.},
\label{wk}
\end{equation}
where $O$ is uniform in $k\in[1,N]$. Moreover,
\begin{equation}
W^*(2\nu_N)-W^*(N)=O(N^{1/4+\varepsilon/2})\quad{\rm a.s.}
\label{w*}
\end{equation}
\begin{equation}
\xi(0,\rho_N)-\xi(0,\rho_{\nu_N}^+)=O(\log N)\quad{\rm a.s.}
\label{xirho}
\end{equation}
\end{lemma}

\noindent{\bf Proof.}
\begin{eqnarray*}
&&P\left(\bigcup_{k=1}^N  \left\{|W(k,2\nu_N)-W(k,N)|\geq
k^{1/2}N^{1/4+\varepsilon/2 }\right\}\right)\\
&&\leq\sum_{k=1}^N P\left(\frac{|W(k,2\nu_N)-W(k,N)|}{k^{1/2}}\geq
N^{1/4+\varepsilon/2}\right)\\
&&\leq NP\left(\widetilde W(|2\nu_N-N|)\geq
N^{1/4+\varepsilon/2}\right)\\
&&\leq
NP\left(\sup_{0\leq u\leq N^{1/2+\varepsilon/2}}|\widetilde W(u)|\geq
N^{1/4+\varepsilon/2}\right)+NP(|2\nu_N-N|\geq N^{1/2+\varepsilon/2})\\
&&\leq 2N\exp(-N^{\varepsilon/2}/2)+2N\exp(-N^\varepsilon/2),
\end{eqnarray*}
where $\widetilde W(\cdot)$ is a standard Wiener process and we used
Lemmas 3.3 and 3.6 (with $k=1$).

Hence (\ref{wk}) follows by Borel-Cantelli lemma, and
(\ref{w*}) follows from (\ref{wk}) by putting $k=1$ there.
To show (\ref{xirho}), observe that
$$
P(\xi(0,\rho_N)-\xi(0, \rho_{\nu_N}^+)\geq j)=\frac{1}{2^j},
$$
since the event $\{\xi(0,\rho_N)-\xi(0, \rho_{\nu_N}^+)\geq j\}$ means
that the last $j$ excursions out of $N$ are downward and, looking at the
random walk from $\rho_N$ backward, this event is equivalent to the
event that the first $j$ excursions are downward, which has the
probability $1/2^j$. Putting $j=2\log N$, (\ref{xirho}) follows by
Borel-Cantelli lemma.
$\Box$

This also completes the proof of Proposition 1.3.
$\Box$

\medskip
Now we are ready to prove Theorem 1.1. Put $N=\xi(0,n)$ into
(\ref{inv7}).  By Proposition 1.3 we have
\begin{equation}
\xi(k,\kappa_n)-\xi(0,n)=G(k,\xi(0,n))+
O(k^{5/4}(\xi(0,n))^{1/4+\varepsilon/2})\quad{\rm a.s.},
\label{xikappa}
\end{equation}
where $\kappa_n=\max\{i\leq n:\, S_i=0\}$, i.e.,
the last zero before $n$ of the random walk and $O$ is uniform for
$k\in [1,(\xi(0,n))^{1/3-\varepsilon})$.

\begin{lemma} For any $\delta>0$, as $n\to\infty$,
\begin{equation}
\xi(k,n)-\xi(k,\kappa_n)=O(kn^\delta)\quad{\rm a.s.},
\label{xikappa2}
\end{equation}
where $O$ is uniform in $k\in [1,n]$.
\end{lemma}

\noindent{\bf Proof.} We have
$$
\xi(k,n)-\xi(k,\kappa_n)
\leq \max_{0\leq i\leq \xi(0,n)}
(\xi(k,\rho_{i+1})-\xi(k,\rho_i)),
$$
therefore
\begin{eqnarray*}
&&P\left(\bigcup_{k=1}^{n}\{\xi(k,n)-\xi(k,\kappa_n)\geq
kn^\delta\}\right)\\
&&\leq\sum_{k=1}^nP\left(\max_{0\leq i\leq\xi(0,n)}
(\xi(k,\rho_{i+1})-\xi(k,\rho_i))\geq kn^\delta\right)\\
&&\leq P(\xi(0,n)\geq n^{1/2+\delta})
+\sum_{k=1}^n P\left(\max_{0\leq i\leq n^{1/2+\delta}}
(\xi(k,\rho_{i+1})-\xi(k,\rho_i))\geq kn^\delta\right)\\
&&\leq P(\xi(0,n)\geq n^{1/2+\delta})
+n^{1/2+\delta}\sum_{k=1}^n P(\xi(k,\rho_1)\geq kn^\delta).
\end{eqnarray*}

Lemma 3.9 implies
\begin{equation}
P(\xi(0,n)\geq n^{1/2+\delta})\leq C_1e^{-C_2n^{2\delta}}.
\label{xitail}
\end{equation}
Moreover, from the distribution of $\xi(k,\rho_1)$ (cf. R\'ev\'esz
\cite{re2} pp. 100, Theorem 9.7), we get
\begin{equation}
P(\xi(k,\rho_1)\geq j)=\frac1{2k}\left(1-\frac1{2k}\right)^{j-1}
\leq e^{-j/(2k)}.
\label{xikrho}
\end{equation}
Putting $j=kn^\delta$, (\ref{xikappa2}) follows from
(\ref{xitail}) and (\ref{xikrho}) by applying Borel-Cantelli lemma.
$\Box$

To complete the proof of Theorem 1.1, observe that for any $\delta>0$,
almost surely
$$
n^{1/2-\delta}\leq \xi(0,n)\leq n^{1/2+\delta}
$$
for all $n$ large enough. We have, as $n\to\infty$,
$$
(\xi(0,n))^{1/4+\varepsilon/2}=O(n^{1/8+5\varepsilon/8})\quad\mathrm{a.s.}
$$
Now (\ref{gkn}) follows from (\ref{xikappa}) and Lemma 4.4, since for
large $n$ the $O$ term in (\ref{xikappa}) is uniform in
$k\in[1,n^{1/6-\varepsilon})$, as stated.
$\Box$

\section{Proof of Theorems 1.2 and 1.3}
\renewcommand{\thesection}{\arabic{section}} \setcounter{equation}{0}
\setcounter{thm}{0} \setcounter{lemma}{0}

In this section we show that the local time $\xi(0,n)$ in
(\ref{gkn}) can be changed to another random walk local time
$\widetilde\xi(0,n)$ and also to a Brownian local time $\eta(0,n)$, both
independent of $G(\cdot,\cdot)$, as claimed in Theorems 1.2 and 1.3,
respectively. The method of proof is similar to that of \cite{cscsfr1},
\cite{cscsfr2}.

\subsection{Proof of Theorem 1.2}
Assume that on the same probability space we have two independent simple
symmetric random walks $\{S^{(1)}_i,\, i=1,2,\ldots\}$ and
$\{S^{(2)}_i,\, i=1,2,\ldots\}$, with respective local times
$\xi^{(1)}(\cdot,\cdot)$ and $\xi^{(2)}(\cdot,\cdot)$. Assume
furthermore that the above procedure has been performed for both random
walks, i.e. we have Wiener sheets $W^{(1)}(\cdot,\cdot)$,
$W^{(2)}(\cdot,\cdot)$ and Wiener processes $W^{*(1)}$, $W^{*(2)}$
satisfying Propositions 1.1-1.3 and Theorem 1.1. Based on these two
random walks, we construct a new simple symmetric random walk
$\{S_i,\, i=1,2,\ldots\}$ such that its local time $\xi(0,n)$ will be
close to $\xi^{(1)}(0,n)$, while the increments $\xi(k,n)-\xi(0,n)$ will
be close to $\xi^{(2)}(k,n)-\xi^{(2)}(0,n)$. This is achieved by taking
"large" excursions from $S^{(1)}$ and "small" excursions from $S^{(2)}$.
As a result, we shall conclude that $\xi(k,n)-\xi(0,n)$ can be
approximated by $G^{(2)}(k,\xi^{(1)}(0,n))$.

This is done as follows (see \cite{cscsfr2}).
Let $\rho^{(j)}_i,\, j=1,2,\, i=1,2,\ldots$ denote the consecutive
return times to zero of the random walk $S^{(j)}$. Let furthermore
$N_0=0$, $N_\ell=2^\ell$, $r_\ell=N_\ell-N_{\ell-1}=2^{\ell-1}$,
$\ell =1,2,\ldots$, and consider the blocks out of which the $\ell$-th
block consisting of $r_\ell$ excursions as follows.
$$
\left\{S^{(j)}_{\rho^{(j)}_{N_{\ell-1}}+1},\ldots,
S^{(j)}_{\rho^{(j)}_{N_\ell}}\right\},\quad j=1,2,\quad \ell=1,2,\ldots
$$
In this block call an excursion large if
$$
\rho^{(j)}_{N_{\ell-1}+i}-\rho^{(j)}_{N_{\ell-1}+i-1}>r_\ell^{4/3},
$$
and call it small otherwise. Now construct the block
$$
\{S_{\rho_{N_{\ell-1}}+1},\ldots, S_{\rho_{N\ell}}\}
$$
of the new random walk, the $\ell$-th block having also $r_\ell$
excursions by keeping large excursions in the block of $S^{(1)}$ unaltered
and replacing small excursions of $S^{(1)}$ by the small excursions of
$S^{(2)}$, keeping also the order of small and large excursions as it was in
$S^{(1)}$. It is possible that there are more small excursions in the block
of $S^{(1)}$ than in the block of $S^{(2)}$. In this case replace as
many small excursions as possible by those of $S^{(2)}$, and leave the
other small excursions unaltered in $S^{(1)}$. One can easily see that,
putting these blocks one after the other, the resulting $S_1,S_2,\ldots$
is a simple symmetric random walk. We denote by $\xi,\rho$, etc.,
without superfix, the corresponding quantities defined for this random
walk, and continue with establishing the next five lemmas that will also
lead to concluding Theorem 1.2.

\begin{lemma}The following inequalities hold:
\begin{eqnarray}
&&\qquad\max_{1\leq i\leq N_\ell}|\rho_i-\rho^{(1)}_i|
\label{rho1}\\
&&\leq\sum_{j=1}^2\sum_{m=1}^\ell\sum_{i=1}^{r_m}
\left(\rho_{N_{m-1}+i}^{(j)}-\rho_{N_{m-1}+i-1}^{(j)}\right)
I\left\{\rho_{N_{m-1}+i}^{(j)}-\rho_{N_{m-1}+i-1}^{(j)}\leq
r_m^{4/3}\right\},\nonumber
\end{eqnarray}
and
\begin{equation}
\max_{1\leq i\leq N_\ell}|\xi(k,\rho_i)-\xi^{(2)}(k,\rho_i^{(2)})|
\leq \xi^*(k)\sum_{m=1}^\ell (\mu_m^{(1)}+\mu_m^{(2)}),
\label{xi4}
\end{equation}
where $I\{\cdot\}$ denotes the indicator of the event in the brackets,
\begin{equation}
\xi^*(k)
=\max_{j=1,2}\max_{1\leq i\leq N_\ell}
\left(\xi^{(j)}(k,\rho_i^{(j)})-\xi^{(j)}(k,\rho_{i-1}^{(j)})\right)
\end{equation}
and $\mu_m^{(j)}$ is the number of large excursions in the $m$-th block
of $S^{(j)}$.
\end{lemma}

\noindent{\bf Proof.}
Obviously, $\max_{1\leq i\leq N_\ell}|\rho_i-\rho_i^{(1)}|$ can be
overestimated by the total length of small excursions of the two random
walks up to time $N_\ell$ which is the right-hand side of (\ref{rho1}).

Moreover, $|\xi(k,\rho_i)-\xi^{(2)}(k,\rho_i^{(2)})|$ can be
overestimated by the total number of large excursions up to $N_\ell$
multiplied by the maximum of the local time of $k$ over all excursions
up to $N_\ell$ of the two random walks, which is the right-hand side of
(\ref{xi4}).
$\Box$

\begin{lemma}
For $n\leq \rho_N^{(1)}$ we have
\begin{equation}
\max_{1\leq i\leq n}|\xi(0,i)-\xi^{(1)}(0,i)|
\leq \max_{1\leq j\leq N}|\xi(0,\rho_j)-\xi(0,\rho_j^{(1)})|+1.
\end{equation}
\end{lemma}

\noindent{\bf Proof.}
Since $\xi(0,\rho_j)=\xi^{(1)}(0,\rho_j^{(1)})=j$, we have for
$\rho_{j-1}^{(1)}\leq i<\rho_j^{(1)}$, $j\leq N$,
$$
\xi(0,i)-\xi^{(1)}(0,i)\leq \xi(0,\rho_j^{(1)})-(j-1)
$$
$$
=\xi(0,\rho_j^{(1)})-\xi(0,\rho_j)+1\leq
\max_{1\leq j\leq N}|\xi(0,\rho_j)-\xi(0,\rho_j^{(1)})|+1.
$$
On the other hand,
$$
\xi^{(1)}(0,i)-\xi(0,i)\leq j-1-\xi(0,\rho_{j-1}^{(1)})
$$
$$
=\xi(0,\rho_{j-1})-\xi(0,\rho_{j-1}^{(1)})\leq
\max_{1\leq j\leq N}|\xi(0,\rho_j)-\xi(0,\rho_j^{(1)})|+1.
$$
$\Box$

\begin{lemma} For $C>0$, $K=1,2,\ldots$ we have
\begin{eqnarray}
&&P\left(\bigcup_{k=1}^K\{\max_{1\leq i\leq N_\ell}
|\xi(k,\rho_i)-\xi^{(2)}(k,\rho_i^{(2)})|\geq 3Ck\ell^2
r_\ell^{1/3}\}\right)\nonumber\\
&&\leq N_\ell\sum_{k=1}^K\frac1k\left(1-\frac1{2k}\right)^{Ck\log
N_\ell}+K\exp(2(e-3)\ell r_\ell^{1/3}). \label{nh}
\end{eqnarray}
\end{lemma}

\noindent{\bf Proof.} Using (\ref{xi4}) of Lemma 5.1, and $4\log 2<3$,
we get
$$
P\left(\bigcup_{k=1}^K\left\{\max_{1\leq i\leq N_\ell}
|\xi(k,\rho_i)-\xi^{(2)}(k,\rho_i^{(2)})|\geq 3Ck\ell^2
r_\ell^{1/3}\right\}\right)
$$
$$
\leq \sum_{k=1}^K P(\xi^*(k)\geq Ck\log N_\ell)
+KP\left(\sum_{m=1}^\ell (\mu_m^{(1)}+\mu_m^{(2)})\geq 4\ell
r_\ell^{1/3}\right).
$$
Using again the distribution of $\xi(k,\rho_1)$ in \cite{re2}, we get
$$
P(\xi^*(k)\geq Ck\log N_\ell)\leq 2N_\ell P(\xi(k,\rho_1)\geq Ck\log
N_\ell)
$$
$$
\leq\frac{N_\ell}{k}\left(1-\frac1{2k}\right)^{Ck\log N_\ell}.
$$

Moreover, $\{\mu_m^{(j)},\, j=1,2,\, m=1,2,\ldots\}$ are
independent random variables such that $\mu_m^{(1)}+\mu_m^{(2)}$ has
binomial distribution with parameters $2r_m$ and
$p_m=P(\rho_1\geq r_m^{4/3})\leq r_m^{-2/3}$, where Lemma 3.10 was used
for $N=1$. Using the moment generating function of the binomial
distribution and exponential Markov's inequality, proceeding as in
\cite{cscsfr2}, we get
$$
P\left(\sum_{m=1}^\ell (\mu_m^{(1)}+\mu_m^{(2)})\geq z\right)\leq
e^{-z}\prod_{m=1}^\ell (1+p_m(e-1))^{2r_m}
$$
$$
\leq \exp\left(2(e-1)\sum_{m=1}^\ell r_mp_m-z\right)
\leq\exp(2(e-1)\ell r_\ell^{1/3}-z).
$$
Putting $z=4\ell r_\ell^{1/3}$, we get (\ref{nh}).
$\Box$

\begin{lemma}As $N\to\infty$,
\begin{equation}
\xi(k,\rho_N)-\xi^{(2)}(k,\rho_N^{(2)})=O(kN^{1/3}\log^2 N)\quad
\mathrm{a.s.},
\label{xikrho2}
\end{equation}
where $O$ is uniform in $k\in [1,N]$.
\end{lemma}

\noindent{\bf Proof.}
Applying the inequality (\ref{nh}) in Lemma 5.3 with  $K=N_\ell$, the
right hand side is summable for $\ell$, provided that  $C$ is large
enough. Hence
$$
\max_{1\leq i\leq N_\ell}|\xi(k,\rho_i)-\xi^{(2)}(k,\rho_i^{(2)})|
=O(k\ell^2 r_\ell^{1/3})=O(k(\log N_\ell)^2 N_\ell^{1/3})
$$
almost surely, as $\ell\to\infty$, from which (\ref{xikrho2}) follows.
$\Box$

To verify Theorem 1.2, we start from (\ref{inv7}) in Proposition 1.3,
applying it for the random walk $S^{(2)}$. We have
$$
\xi^{(2)}(k,\rho_N^{(2)})-\xi^{(2)}(0,\rho_N^{(2)})=
G^{(2)}(k,N)+O(k^{5/4}N^{1/4+\varepsilon/2})\quad\mathrm{a.s.}
$$
as $N\to\infty$. Since $\xi^{(2)}(0,\rho_N^{(2)})=\xi(0,\rho_N)=N$,
according to Lemma 5.4 we also have, as $N\to\infty$,
$$
\xi(k,\rho_N)-\xi(0,\rho_N)=
G^{(2)}(k,N)+O(k^{5/4}N^{1/4+\varepsilon/2}+kN^{1/3}\log^2 N)
$$
almost surely. Now put $N=\xi(0,n)$. Using Lemma 4.4, we can see as
before,
$$
\xi(k,n)-\xi(0,n)=G^{(2)}(k,\xi(0,n))+
O(k^{5/4}n^{1/8+5\varepsilon/8}+kn^{1/6+\varepsilon/4})
$$
almost surely and uniformly in $k\in [1,n^{1/6-\varepsilon}]$, as
$n\to\infty$. It remains to show that on the
right-hand side $\xi(0,n)$ can be replaced by $\xi^{(1)}(0,n)$.

\begin{lemma} For any $\varepsilon>0$ there exists a $\delta>0$ such
that, as $n\to\infty$,
\begin{equation}
|G^{(2)}(k,\xi(0,n))-G^{(2)}(k,\xi^{(1)}(0,n))|=O(k^{1/2}n^{1/4-\delta})
\quad {\rm a.s.},
\label{gxi1}
\end{equation}
where $O$ is uniform in $k\in[1,n^{1/6-\varepsilon}]$.
\end{lemma}

\noindent{\bf Proof.}
Let $0<\varepsilon<1/6$ and $K_\ell=[2^{\ell(1/6-\varepsilon)}]$,
$u_\ell=2^{\ell(1/4-\varepsilon/100)}$. Since $k^{-1/2}W(k,\cdot)$
is a standard Wiener process (denoted by $\widetilde W(\cdot)$),
we have
\begin{eqnarray*}
&&P\left(\bigcup_{k=1}^{K_\ell} \left\{\max_{2^{\ell-1}\leq
n<2^\ell} |W^{(2)}(k,\xi(0,n))-W^{(2)}(k,\xi^{(1)}(0,n))| \geq
k^{1/2}u_\ell\right\}\right)\\
&&\leq K_\ell P\left(\max_{2^{\ell-1}\leq n<2^\ell} |\widetilde
W(\xi(0,n))-\widetilde W(\xi^{(1)}(0,n)| \geq u_\ell\right)\\
&&\leq K_\ell P(\sup_{(u,v)\in A}|\widetilde W(u)-\widetilde W(v)|
\geq u_\ell)\\
&&+2K_\ell P(\xi(0,2^\ell)\geq 2^\ell)\\
&&+K_\ell P\left(\max_{1\leq n\leq 2^\ell}
|\xi(0,n)-\xi^{(1)}(0,n)|\geq 2^{\ell(1/2-\varepsilon/48)}\right),
\end{eqnarray*}
where
$$
A=\{(u,v):\, 0\leq u\leq 2^\ell,\, 0\leq v\leq 2^\ell, |u-v|\leq
2^{\ell(1/2-\varepsilon/48)}\}.
$$

First we estimate the last term. By Lemma 5.2
\begin{eqnarray*}
&&P\left(\max_{1\leq n\leq 2^\ell} |\xi(0,n)-\xi^{(1)}(0,n)|\geq
2^{\ell(1/2-\varepsilon/48)}\right)\\
&&\leq P\left(\max_{1\leq j\leq 2^{\ell(1/2+\varepsilon/4)}}
|\xi(0,\rho_j)-\xi(0,\rho^{(1)}_j)|\geq
2^{\ell(1/2-\varepsilon/48)}-1\right)\\
&&+P\left(\rho^{(1)}_{[2^{\ell(1/2+\varepsilon/4)}]}\leq 2^\ell\right)\\
&&
\leq P\left(\max_{(i,j)\in B}|\xi(0,i)-\xi(0,j)|\geq
2^{\ell(1/2-\varepsilon/48)}-1\right)\\
&&
+P\left(\rho^{(1)}_{[2^{\ell(1/2+\varepsilon/4)}]}\leq 2^\ell\right)
+2P\left(\rho_{[2^{\ell(1/2+\varepsilon/4)}]}\geq
2^{\ell(4/3+\varepsilon)}\right)\\
&&
+P\left(\max_{1\leq j\leq
2^{\ell(1/2+\varepsilon/4)}}|\rho_j-\rho_j^{(1)}| \geq
2^{\ell(1-\varepsilon/12)}\right),
\end{eqnarray*}
where
$$
B=\{(i,j):\, 1\leq i\leq 2^{\ell(4/3+\varepsilon)}, 1\leq j\leq
2^{\ell(4/3+\varepsilon)}, |i-j|\leq 2^{\ell(1-\varepsilon/12)}\}.
$$

Now we estimate the respective right-hand sides of the previous
two inequalities term by term.

Lemma 3.9 implies
$$
P(\sup_{(u,v)\in A}|\widetilde W(u)-\widetilde W(v)| \geq
2^{\ell(1/4-\varepsilon/100)})
$$
$$
\leq C_12^{\ell(1/2+\varepsilon/48)}\exp\left(-C_2
2^{\ell(\varepsilon/48-\varepsilon/50)}\right),
$$
and
$$
P\left(\max_{(i,j)\in B}|\xi(0,i)-\xi(0,j)|\geq
2^{\ell(1/2-\varepsilon/48)}-1\right)
$$
$$
\leq C_1
2^{\ell(1/6+13\varepsilon/24)}\exp\left(-C_2
(2^{\ell\varepsilon/24}-2)\right).
$$
Observe that
$$
P(\xi(0,2^\ell)\geq 2^\ell)=0
$$
and
$$
P(\rho^{(1)}_{[2^{\ell(1/2+\varepsilon/4)}]}\leq 2^\ell)=
P(\xi^{(1)}(0,2^\ell)\geq 2^{\ell(1/2+\varepsilon/4)}) \leq C_1e^{-C_2
2^{\ell\varepsilon/2}}.
$$

From Lemma 3.10 we have
$$
P(\rho^{(1)}_{[2^{\ell(1/2+\varepsilon/4)}]}\geq
2^{\ell(4/3+\varepsilon)}) \leq C2^{-\ell(1/6+\varepsilon/4)}.
$$

Finally, from (\ref{rho1}) of Lemma 5.1, Lemma 3.10 and Markov's
inequality
$$
P\left(\max_{1\leq j\leq
2^{\ell(1/2+\varepsilon/4)}}|\rho_j-\rho_j^{(1)}| \geq
2^{\ell(1-\varepsilon/12)}\right)
$$
$$
\leq
\frac{2}{2^{\ell(1-\varepsilon/12)}}\sum_{m=1}^{\ell(1/2+\varepsilon/4)}
r_m E(\rho_1 I(\rho_1\leq r_m^{4/3}))
$$
$$
\leq \frac{C}{2^{\ell(1-\varepsilon/12)}}
\sum_{m=1}^{\ell(1/2+\varepsilon/4)}2^{5(m-1)/3}\leq
C2^{\ell(-1/6+\varepsilon/2)}.
$$

Assembling all these estimations, we obtain
$$
P\left(\bigcup_{k=1}^{K_\ell} \left\{\max_{2^{\ell-1}\leq
n<2^\ell} |W^{(2)}(k,\xi(0,n))-W^{(2)}(k,\xi^{(1)}(0,n))| \geq
k^{1/2}u_\ell\right\}\right)
$$
$$
\leq C_12^{2\ell /3}
\exp\left(-C_22^{\ell\varepsilon(1/48-1/50)}\right)+ C_3
2^{\ell/6}\exp(-C_22^{\ell\varepsilon/2})
$$
$$
+C_12^{\ell(1/3-11\varepsilon/24)}
\exp\left(-C_2(2^{\ell\varepsilon/24}-2)\right)
+C2^{-5\ell \varepsilon/4}.
$$

Since all these terms are summable in $\ell$, by Borel-Cantelli lemma we
have
$$
\max_{2^{\ell-1}\leq n<2^\ell}
|W^{(2)}(k,\xi(0,n))-W^{(2)}(k,\xi^{(1)}(0,n))|
=O(k^{1/2}2^{\ell(1/4-\varepsilon/100)})
$$
almost surely, as $\ell\to\infty$, uniformly for $k\in
[1,2^{\ell(1/6-\varepsilon)}]$, i.e.,
$$
|W^{(2)}(k,\xi(0,n))-W^{(2)}(k,\xi^{(1)}(0,n))|
=O(k^{1/2}n^{1/4-\varepsilon/100})
$$
almost surely, as $n\to\infty$, uniformly for $k\in
[1,n^{1/6-\varepsilon}]$.
Similar estimations hold for the other terms of $G^{(2)}$, hence
we have (\ref{gxi1}) with $\delta=\varepsilon/100$.
$\Box$

Since the above estimations also imply
$$
\xi(0,n)-\xi^{(1)}(0,n)=O(n^{1/2-\delta})
$$
almost surely, when $n\to\infty$, with $\delta=\varepsilon/48$,
on choosing $\widetilde\xi(0,\cdot)=\xi^{(1)}(0,\cdot)$,
$G(\cdot,\cdot)=G^{(2)}(\cdot,\cdot)$, the proof of
Theorem 1.2 is completed as well.
$\Box$

\subsection{Proof of Theorem 1.3}

First, we give a coupling inequality for the invariance principle
between random walk and Brownian local times at location zero. We use
Skorokhod embedding as in \cite{csh}, i.e., given a standard Wiener
process $W(\cdot)$ with its local time $\eta(\cdot,\cdot)$, define a
sequence of stopping times $\{\tau_i\}_{i=0}^\infty$ by $\tau_0=0$,
$$
\tau_n:=\inf\{t:\, t>\tau_{n-1},\, |W(t)-W(\tau_{n-1})|=1\},
\quad n=1,2,\ldots
$$
Then $S_n=W(\tau_n),\, n=0,1,2,\ldots$ is a simple symmetric random
walk. Denote by $\xi(\cdot,\cdot)$ its local time and by $\rho_i$ the
return times to zero. Moreover, define
$$
\eta_i:=\eta(0,\tau_{\rho_i+1})-\eta(0,\tau_{\rho_i}),
$$
i.e., the Brownian local time between the $i$-th return to zero and
next stopping time $\tau$. Then by Knight \cite{kni} the random
variables $\eta_i,\, i=1,2,$ are i.i.d. having exponential
distribution with parameter $1$. There is no other contribution than
$\eta_i$ to the Brownian local time $\eta(0,\cdot)$. Moreover, we have
$$
\left|\eta(0,\tau_n)-\sum_{i=1}^{\xi(0,n)}\eta_i\right|\leq
\eta_{\xi(0,n)},
$$
the error term being zero if $S_n=W(\tau_n)\neq 0$. If $S_n=0$, then the
last term $\eta_{\xi(0,n)}$ is not counted in $\eta(0,\tau_n)$.
Now we have
$$
|\xi(0,n)-\eta(0,n)|\leq |\eta(0,\tau_n)-\eta(0,n)|
+\max_{1\leq j\leq \xi(0,n)}\left|\sum_{i=1}^j
(\eta_i-1)\right|+\eta_{\xi(0,n)}.
$$
Therefore, for $\delta>0$
\begin{eqnarray*}
&&P(|\xi(0,n)-\eta(0,n)|\geq 2n^{1/4+\delta}+C\log n)\\
&&\leq P(\xi(0,n)\geq n^{1/2+\delta})
+P\left(\max_{1\leq j\leq
n^{1/2+\delta}}\left|\sum_{i=1}^j(\eta_i-1)\right|\geq
n^{1/4+\delta}\right)\\
&&+P\left(\max_{1\leq i\leq n^{1/2+\delta}}\eta_i\geq C\log n\right)
+P(|\tau_n-n|\geq n^{1/2+\delta})\\
&&+P\left(\sup_{|u-n|\leq n^{1/2+\delta}}|\eta(0,u)-\eta(0,n)|\geq
n^{1/4+\delta}\right).
\end{eqnarray*}
Estimating the above probabilities term by term, by Lemmas 3.9, 3.11 and 3.12,
\begin{eqnarray*}
&&P(\xi(0,n)\geq n^{1/2+\delta})\leq C_1e^{-C_2 n^{2\delta}},\\
&&P\left(\max_{1\leq j\leq
n^{1/2+\delta}}\left|\sum_{i=1}^j(\eta_i-1)\right|\geq
n^{1/4+\delta}\right)\leq 2e^{-n^\delta/8},\\
&&P\left(\max_{1\leq i\leq n^{1/2+\delta}}\eta_i\geq C\log n\right)\leq
n^{1/2+\delta-C},\\
&&P(|\tau_n-n|\geq n^{1/2+\delta})\leq 2e^{-3n^{2\delta}/8},\\
&&P\left(\sup_{|u-n|\leq n^{1/2+\delta}}|\eta(0,u)-\eta(0,n)|\geq
n^{1/4+\delta}\right)\leq C_1n^{1/4-\delta/2}e^{-C_2n^{\delta}}.
\end{eqnarray*}

Hence, we arrive at the coupling inequality for the invariance principle
between random walk and Brownian local times
\begin{eqnarray}
&&P(|\xi(0,n)-\eta(0,n)|\geq 2n^{1/4+\delta}+C\log n)\nonumber\\
&&\leq C_1n^{1/4-\delta/2}e^{-C_2n^{\delta}}+n^{1/2+\delta-C}.
\label{xieta}
\end{eqnarray}

By choosing $0<\delta<1/2$ and $C>2$, this also implies
\begin{equation}
\xi(0,n)-\eta(0,n)=O(n^{1/4+\delta})\quad\mathrm{a.s.}
\label{invxieta}
\end{equation}
as $n\to\infty$.

For the proof of Theorem 1.3 we apply the above procedure for
$\xi^{(1)}(0,\cdot)$, i.e., we construct a standard Brownian local time
$\eta(0,\cdot)$ satisfying the above inequality with $\xi$ replaced by
$\xi^{(1)}$. We may assume that $\eta(0,\cdot)$ is also independent of
$G(\cdot,\cdot)$ of Theorem 1.2. We show that in (iv) of Theorem 1.2,
$\widetilde\xi=\xi^{(1)}$ can be replaced by $\eta$ with the same $O$
term.

\begin{lemma} As $n\to\infty$, we have for any $\delta>0$
$$
\left|W^{(2)}(k, \xi^{(1)}(0,n))-
W^{(2)}(k,\eta(0,n))\right|=O(k^{1/2} n^{1/8+\delta})\quad\mathrm{a.s.},
$$
where $O$ is uniform in $k\in[1,n^{1/6}]$.
\end{lemma}

\noindent{\bf Proof.} Let $K_n=[n^{1/6}]$.
\begin{eqnarray*}
&&P\left(\bigcup_{k=1}^{K_n} \left|W^{(2)}(k, \xi^{(1)}(0,n))-
W^{(2)}(k,\eta(0,n)))\right|\geq k^{1/2} n^{1/8+\delta}\right)
\\
&&\leq K_n P(\sup_{(u,v)\in D}|\widetilde W(u)-\widetilde W(v)|\geq
n^{1/8+\delta}) +K_n P(\xi^{(1)}(0,n)\geq n^{1/2+\delta})\\
&&+K_n P(\eta(0,n)\geq n^{1/2+\delta})\\
&&+K_n P(|\xi^{(1)}(0,n)-\eta(0,n)|\geq 2n^{1/4+\delta}+C\log n),
\end{eqnarray*}
where
$$
D=\{(u,v):\, u\leq n^{1/2+\delta},\, v\leq n^{1/2+\delta},
\, |u-v|\leq 2n^{1/4+\delta}+C\log n\}.
$$
Now using Lemma 3.9, we get the inequalities
$$
K_nP\left(\sup_{(u,v)\in D}|\widetilde W(u)-\widetilde W(v)|\geq
n^{1/8+\delta}\right)\leq C_1 n^{1/2}e^{-C_2n^\delta},
$$
$$
K_n P(\xi^{(1)}(0,n)\geq n^{1/2+\delta})\leq
C_1n^{1/6}e^{-C_2n^{2\delta}},
$$
$$
K_n P(\eta(0,n)\geq n^{1/2+\delta})\leq C_1 n^{1/6}e^{-C_2n^{2\delta}}.
$$

By choosing $C$ large enough in (\ref{xieta}), the right-hand sides of
that and also of the above inequalities are summable in $n$, Lemma 5.6
follows by Borel-Cantelli lemma.
$\Box$

The same holds for other terms of $G^{(2)}$. Choosing
$\delta<\varepsilon$, the error term in Lemma 5.6 is smaller than
$kn^{1/6+\varepsilon}$ in (iv) of Theorem 1.2, hence we have also (iii)
of Theorem 1.3 with $G=G^{(2)}$, $\widetilde \xi=\xi^{(1)}$. (ii) of
Theorem 1.3 follows from (\ref{invxieta}) with $\xi$ replaced by
$\xi^{(1)}=\widetilde \xi$ and (iii) of Theorem 1.2. This completes the
proof of Theorem 1.3.
$\Box$

\section*{Acknowledgements}
We are grateful to Zhan Shi for valuable suggestions.

\end{document}